\documentclass[a4paper,11pt,leqno]{article}
\usepackage{authblk}
\usepackage{stmaryrd}
\usepackage[utf8]{inputenc}
\usepackage[numbers,sort&compress]{natbib}
\usepackage{amsmath,amsthm}
\usepackage[normalem]{ulem}
\usepackage{amssymb,latexsym}
\usepackage{soul}

\usepackage{mwe}
\usepackage{caption}
\usepackage{mathrsfs}
\DeclareMathOperator*{\esssup}{ess\,sup}

\newcommand{\R}{\mathbb{R}}

\newcommand{\Q}{\mathbb{Q}}

\newcommand{\dis}{\displaystyle}

\numberwithin{equation}{section}
\usepackage[
left=34mm,
right=34mm,
top=25mm,
bottom=25mm,
includeheadfoot,      
]{geometry}

\usepackage{xcolor}

\theoremstyle{plain}
\newtheorem{theorem}{Theorem}[section]
\newtheorem{proposition}{Proposition}[section]
\newtheorem{definition}{Definition}[section]

\newtheorem{remark}{Remark}[section]
\newtheorem{lemma}{Lemma}[section]

\title{A class of parabolic  reaction-diffusion systems governed by spectral fractional Laplacians : Analysis and numerical simulations}
\author{Maha DAOUD}
\affil{ Département Génie Mathématique et Modélisation, INSA Toulouse, 31400 Toulouse, France\\\vspace*{.15cm} mahaadaoud@gmail.com}
\newcommand{\midsize}{\fontsize{9.35pt}{11.5pt}\selectfont}
\begin{document}

\maketitle

\begin{abstract}
In this paper, we prove the global-in-time existence of strong solutions to a class of fractional parabolic reaction-diffusion systems set in a bounded open subset of $\mathbb{R}^N$. The diffusion operators are of the form $u_i \mapsto d_i (-\Delta)_{Sp}^{s_i} u_i$, where $0 < s_i < 1$. The operator $(-\Delta)_{Sp}^{s}$ stands for the commonly called  spectral fractional Laplacian. Moreover, the nonlinear reaction terms are assumed to fulfill natural structural conditions that ensure the nonnegativity of the solutions and provide uniform control of the total mass.
We establish the global existence of strong solutions under the assumption that the nonlinearities exhibit at most polynomial growth. 		 
Our results extend previous results obtained when the diffusion operators are of the form \( u_i \mapsto d_i (-\Delta)^s u_i \), where $(-\Delta)^s$ denotes the widely known regional fractional Laplacian.		 
Furthermore, we present some numerical simulations to address a theoretical question that remains open to date.		

\bigbreak
\noindent
\textbf{Keywords.} {Reaction-diffusion system, fractional diffusion, strong solution, global existence, numerical simulation.
}

\bigbreak
\noindent
	\textbf{Mathematics Subject Classification (2020):} {35R11, 	35A01,	35D35, 47D06, 35B45.}
\end{abstract}

\section{Introduction}

\maketitle
\label{intro}
 \raggedbottom
In this work, we study the global-in-time existence of nonnegative strong solutions to a class of fractional parabolic reaction-diffusion systems. These systems are governed by spectral fractional Laplacians of different orders and take the following form:
 {\midsize \begin{equation}\label{ReactiondiffusionSystem1}\tag{S}
		\left\{
		\begin{array}{rcll}
			\multicolumn{4}{l}{\forall i=1,\ldots,m,} \\[2pt]
			\partial_t u_i(t,\mathbf{x})+d_i(-\Delta)_{Sp}^{s_i}u_i(t,\mathbf{x})
			&=& f_i\bigl(u_1(t,\mathbf{x}),\ldots,u_m(t,\mathbf{x})\bigr),
			& (t,\mathbf{x})\in(0,T)\times\Omega,\\
			\mathcal{B}[u_i(t,\mathbf{x})] &=& 0,
			& (t,\mathbf{x})\in(0,T)\times\partial\Omega,\\
			u_i(0,\mathbf{x}) &=& u_{0i}(\mathbf{x}),
			& \mathbf{x}\in\Omega,
		\end{array}
		\right.
\end{equation}}
where $\Omega$ is a bounded regular open subset of $\R^N$ with $N\in\mathbb{N}^*$, $m\geq 2$ and for any $i\in\llbracket1, m\rrbracket$, $0<s_i<1$, $d_i>0$ and $f_i$ is locally Lipschitz continuous. {The boundary operator $\mathcal{B}$ stands for either
	$
	\mathcal{B}[u_i]=u_i$ %\quad \text{(homogeneous Dirichlet)}, 
	or
	$\mathcal{B}[u_i]=\frac{\partial u_i}{\partial\nu}$,
	where $\nu$ denotes the outward unit normal to $\partial\Omega$.
	In both cases, $(-\Delta)_{Sp}^{s_i}$ denotes the spectral fractional Laplacian
	associated with the chosen boundary condition.}

\smallbreak \noindent
Let $\{\lambda_k,e_k\}_{k\in\mathbb{N}^*}$ be  the eigenpairs of the ``classical'' Dirichlet or Neumann Laplacian. Then, the corresponding spectral fractional Laplacian is defined for $s\in(0,1)$ by
\begin{eqnarray}\label{Definitionof-}
	(-\Delta)_{Sp}^s u(\mathbf{x}):=\sum\limits_{k=1}^{+\infty} \lambda_{k}^su_{k}e_{k} (\mathbf{x}),
\end{eqnarray}
where $u_k=\displaystyle\int_\Omega u(\mathbf{x}) e_k (\mathbf{x}) d\mathbf{x}$ and $\|e_k\|_{L^2(\Omega)}=1$.

\smallbreak
Parabolic reaction-diffusion systems have long been employed to model a wide range of phenomena, from chemical reactions to biological processes (see, for instance \cite{PierreSurvey2010,SchFis,Pis,Mur} and references therein).  The ``classical'' Laplacian operator, commonly used to describe diffusion, has proven effective in many contexts. Further discussion can be found in \cite{PierreSurvey2010,LaamActa2011,FellLaam2016,LaamPert,LaamPierANHP,GoudVass10,PierSuzYam2019} and their bibliographies. However, in many cases, this operator fails to capture nonlocal interactions that are critical  for accurately simulating complex phenomena.	
This has led to a surge of interest in fractional-type Laplacian operators, which offer a more suitable framework for modeling nonlocal diffusive processes. For further details, we refer the interested reader, for instance, to \cite{SomBur2018,Valdinoci2009,AbaVal2019,DaouLaam,BucVal2016,Vaz2012,Vaz2017,LisAl} and references cited therein. More broadly, fractional partial differential equations have gained prominence in recent years and have become fundamental in many scientific areas, see \cite{Pod,KilSri} and their references.

\smallbreak
In open bounded subsets of $\mathbb{R}^N$, there are various fractional Laplacians,  each with its own characteristics and applications.  Among these different operators, two main classes are widely studied: the spectral fractional Laplacian {(\it SFL)} (typically associated with homogeneous Dirichlet  boundary condition) and the regional fractional Laplacian {\it RFL} defined as 	\begin{equation} \label{FractionalLaplacianDefinition}
	(-\Delta)^su(\mathbf{x}):=a_{N,s}\, \text{\text{P.V.}} \int_{\mathbb{R}^{N}} \frac{u(\mathbf{x})-u(\mathbf{y})}{\|\mathbf{x}-\mathbf{y}\|^{N+2s}}d\mathbf{y},\quad \mathbf{x}\in\Omega,
\end{equation}
where $u(\mathbf{x})=0$ in $\mathbb{R}^N\setminus\Omega$, $a_{N,s}:=\frac{s 2^{2 s} \Gamma\left(\frac{N+2 s}{2}\right)}{\pi^{\frac{N}{2}} \Gamma(1-s)}$, $\|\cdot\|$ is the Euclidean norm of $\R^N$ and \text{P.V.} stands for the Cauchy principal value. While both operators generalize the classical Laplacian, they differ in that one incorporates information from the boundary, while the other considers influences from the exterior of the domain, leading to different stochastic processes. 
For more details about the striking differences between the two operators, see \cite{SerVal} (see also \cite[Section 5.3]{MolRadRaf}, \cite{AbaVal2019,GarSte2019} and references therein). Furthermore, the  author and her collaborator exhibit in \cite[Section 4]{DaouLaam} several simple yet significant one-dimensional examples  that  effectively  highlight the difference between the two operators.  In this context, the same authors and their collaborator have previously studied in \cite{DaouLaamBaal2024} System $(S)$ with the {\it RFL}.
Overall, we have observed that the two operators  lead to notable differences  in the properties of the systems they govern:   they require the adoption of tailored methodologies for each operator and careful consideration of specific technical aspects.

\smallbreak
Recently, considerable studies have explored the existence of solutions to fractional reaction-diffusion problems, particularly for single parabolic equations ($i.e.\;m = 1$) in open bounded subsets of $\mathbb{R}^N$. For further details, we refer the reader to \cite{DipGiaVal2022,ChanGomVz,KobaMats,NaneTuan} for the {\it SFL} and to \cite{GalWarma1,GalWarma2,BicWarZua2,BicWarZua3,FerRos}
for the {\it RFL}.
However, there are relatively few works addressing systems $(m\geq  2)$ governed by fractional Laplacians. For systems with $m=2$ involving the {\it SFL}, we refer to \cite{ZhaoRuan,ShiLiYang}, while for the {\it RFL}, relevant studies can be found in \cite{AtmaBirDaouLaam-Parabo,AtmaBirDaouLaam-Acta,AtmaBirDaouLaam-s1-s2-Ellip,Louis}. As far as we know, the only existing work for $m\geq 2$ is \cite{DaouLaamBaal2024},  where the author and her collaborators analyzed System (S) with the $RFL$ and $s_1=\cdots=s_m$.

	One should question the intriguing case of considering  fractional Laplacians of different orders, $i.e.$ at least for two different $i,j\in\llbracket1,m\rrbracket$ $s_i\neq s_j$.
	This is motivated, for example, by 	stochastic considerations.	
	Specifically, previous studies have investigated equations involving a sum of two fractional Laplacians with different orders. Such sums naturally arise from the superposition of two stochastic processes, each with distinct random walks and Lévy flights. This approach has numerous applications, including chemical reaction design, plasma physics, biophysics and population dynamics. It also serves to model heart anomalies caused by arterial issues.
	For more details, we refer to \cite{ChenBhakHaj2022,DipLiVal2024} and references included therein. For a biological motivation, we refer to \cite{ZhaoRuan}.

\smallbreak	
Let us go back to System (S) with the $SFL$. In this paper, our main goal is to extend  the results obtained in \cite{DaouLaamBaal2024} to a more general setting, specifically  by incorporating different types of  fractional operators.

Since the unknown \(\mathbf{u} = (u_1, \cdots, u_m)\) represents quantities such as chemical concentrations or population densities, we focus on nonnegative solutions. Therefore, the initial data  have to be chosen nonnegative, i.e., \(u_{0i} \geq 0\) for each $i\in\llbracket1,m\rrbracket$. Moreover, it is well-known that solutions, as long as they exist, remain nonnegative if the reaction terms \(f_i\) satisfy the so-called ``quasi-positivity" property, that is:	
{\midsize \begin{equation} \forall i\in\llbracket1,m\rrbracket,\; f_i(r_1,\cdots,r_{i-1},0,r_{i+1},\cdots,r_m) \geq  0,
		\; \forall \mathbf{r}=(r_1,\cdots,r_m) \in [0,+\infty)^m.
		\tag{{\bf P}}
\end{equation}}
In addition, in many cases, the model inherently includes a control or even preservation of the total mass, specifically
\begin{equation}\label{ControlMass}
	\sum_{i=1}^m\int_\Omega u_i(t,\mathbf{x})\,d\mathbf{x} \leq C \sum_{i=1}^m\int_\Omega u_{0i}(\mathbf{x})\,d\mathbf{x} \quad\text{for some } C>0.
\end{equation}
The control (\ref{ControlMass})  is  satisfied if
\begin{equation}
	\exists (a_1,\cdots a_m)\in (0,+\infty)^m \text{ such that } \dis\sum\limits_{i=1}^m a_if_i\leq 0,
	\tag{{\bf M}}
\end{equation}
or, more broadly, if this sum grows at most linearly with respect to its variables, 
{{\it namely}
	\begin{equation}
		\forall \mathbf{r}=(r_1,\ldots,r_m) \in[0,+\infty)^m,\; \sum\limits_{i=1}^m a_i f_i(\mathbf{r}) \leq C\Big[1+\sum\limits_{i=1}^m r_i\Big] \quad\text{for some } C\geq 0.
		\tag{{\bf M'}}
	\end{equation}
}	
\medbreak

It is worth highlighting that properties {\bf(P)} and {\bf(M)} (or {\bf (M')}) are frequently observed in various applications, including models that describe evolutionary processes involving spatial diffusion and chemical reactions. Nevertheless, it turns out that these properties do not keep the solution from blowing up in $L^{\infty}$-norm, even in finite time (see \cite{PierreSchmitt2000}). Hence, in addition to the structure {\bf(P)+ (M)}, some growth restrictions and extra structure on the nonlinearities are needed if one expects global existence of strong solutions.  This issue has been intensively studied in the case of the classical Laplacian, especially when the initial data are bounded. However, many open and challenging problems remain. Indeed, the construction of global solutions is primarily hindered by two major obstructions: (i) when the
the $f_i$'s have quadratic, even faster, growth for large values of the $u_i$'s ; (ii) when the $d_i$'s are very different.

\smallbreak

To situate our work within the existing body of research and underscore the novelty of this paper, let us provide a concise review of the literature concerning the global existence of strong solutions to system \eqref{ReactiondiffusionSystem1}.
\begin{itemize}
	\item  {\it Case where the diffusions are governed by the classical Laplacian (i.e. all $s_i$'s = 1)} 
	
	As mentioned above, for bounded initial data, some sufficient conditions on the fi’s are required to ensure the global existence of  strong solutions (see Definition \ref{Strong_solution}). More precisely,
	\begin{itemize}
		\item[(i)] polynomial growth and the so-called {\it triangular structure} (see~{\eqref{TriangStructure}}). For further details, see~\cite{PierreSurvey2010} and its references;
		\item[(ii)] quadratic growth. For this setting, we refer, $e.g.$, to~\cite{GoudVass10,PierSuzYam2019,Soup18};
		\item[(iii)] when the growth of the~$f_i$'s barely exceed polynomial growth, limited research has been conducted, and solely for $m=2$. The most thoroughly examined model is $f_2(u_1,u_2)=-f_1(u_1,u_2)=u_1e^{u_2^\beta}$. In this case, the global existence of strong solutions has been established in \cite{HaraYouk} for $\beta~<~1$, and in \cite{Barabanova} for $\beta=1$ yet under a restrictive assumption on the size of $u_{01}$. However, the problem remains open for $\beta>1$.
	\end{itemize} 
	Let us emphasize that our primary interest lies in examining strong solutions. It goes without saying that there is a substantial body of research on weak solutions, including notable works such as  \cite{PierreL1,LaamThese1988,Pierre2003,PierreSurvey2010,LaamPierANHP,LaamPert,LaamPierDCDS,LaamPierM3AS}. On the other hand, we would like to underscore that the cited references do not encompass the extensive literature on the subject. For an exhaustive review of findings up to 2010, we refer interested readers  to~\cite{PierreSurvey2010}. This survey not only provides a general overview of the problem but also offers an additional insightful commentary on the mathematical challenges posed by such systems. For some further recent results, see~\cite{LaamActa2011,FellLaam2016,CanDesFel2014,PierSuzYam2019,QuitSoup,Soup18} and references included therein.
	
	\item  {\it Case where the diffusions are governed by fractional Laplacians (i.e. $0<s_i<1$)} 
	\begin{itemize}
		\item[$\diamond$] \textbf{With regional fractional Laplacians}
		
		\begin{itemize}
			
			\item[---] \textit{Systems with two equations ($i.e.\;m = 2$) :}

			In the case where $\Omega$ is bounded, the author and her coworkers have studied  the following System 	
			
			\begin{equation}\label{ABDL}
				\left\{
				\begin{array}{rcll}
					\partial_t u_1 + (-\Delta)^{s_1}u_1 &=& \|\nabla u_2\|^{q}+ h_1 & \text{in }(0,T)\times\Omega,\\
					\partial_t u_2 + (-\Delta)^{s_2}u_2 &=& \|\nabla u_1\|^{p}+ h_2 & \text{in }(0,T)\times\Omega\\
					u_1 = u_2 &=& 0, & \text{in }(0,T)\times(\mathbb{R}^N\setminus\Omega),\\
					u_1(0,\mathbf{x}) &=& u_{01}(\mathbf{x})\ge 0, & \mathbf{x}\in\Omega,\\
					u_2(0,\mathbf{x}) &=& u_{02}(\mathbf{x})\ge 0, & \mathbf{x}\in\Omega.
				\end{array}
				\right.
			\end{equation}
%									
%					\begin{equation} \label{ABDL}
%					\left\{
%					\begin{array}{rclll}
%						\partial_t u_1(t,\mathbf{x})+(-\Delta)^{s_1} u_1(t,\mathbf{x}) &= &\|\nabla u_2(t,\mathbf{x})\|^{q}+ h_1(t,\mathbf{x}), & (t,\mathbf{x})\in (0,T)\times\Omega, \\
%						\partial_t u_2(t,\mathbf{x})+(-\Delta)^{s_2} u_2(t,\mathbf{x}) &= & \|\nabla u_1(t,\mathbf{x})\|^{p}+ h_2(t,\mathbf{x}),  &(t,\mathbf{x})\in (0,T)\times\Omega, \\
%						u_1(t,\mathbf{x}) = u_2(t,\mathbf{x})&=&0, &(t,\mathbf{x})\in(0,T)\times (\R^N\setminus\Omega),\\
%						u_1(0,\textbf{x}) &=& u_{01}(\mathbf{x})\geq 0,  &\mathbf{x}\in\Omega,\\u_2(0,\textbf{x}) &=&u_{02}(\mathbf{x})\geq 0, & \mathbf{x}\in\Omega,
%					\end{array}%
%					\right.
%			\end{equation}	
			where $h_1, h_2\geq 0$ and $p,q\geq1$ (see \cite{AtmaBirDaouLaam-Parabo}). % for the case $s_1=s_2$ and \cite{AtmaBirDaouLaam-s1-s2-Para} for the case $s_1\neq s_2$).						
			Clearly,  the right-hand sides in this case do not fall within the framework of System~\eqref{ReactiondiffusionSystem1}, nor do they satisfy~{\bf(M)}. It is worth highlighting that the same authors have also explored two elliptic versions of System \eqref{ABDL}, namely: the case of identical fractional orders in \cite{AtmaBirDaouLaam-Acta}, and a generalized setting with different orders in \cite{AtmaBirDaouLaam-s1-s2-Ellip}.

			On the other hand, in the case $\Omega=\R^N$, the only works known to us are those of~\cite{Ahmad+4Hindawi,Ahmad+Alsaedi+2,Alsaedi+Al-Yami+2,KiraneAlsaediAhmad}.
			\item[---] \textit{Systems with $m\geq 2$ equations :}
			
			As per our knowledge, the only existing works pertain to the case of equal fractional orders:
			\\ (i) In \cite{DaouLaamBaal2024}, the author and her collaborators investigated System \eqref{ReactiondiffusionSystem1} posed in a bounded open subset $\Omega$. \textit{First,} they have treated the case of reversible chemical reactions of three components (see \eqref{SystemChemReac-DiffFR2-}  below). \textit{Second}, they have examined the existence of strong global nonnegative solutions to System \eqref{ReactiondiffusionSystem1} along with properties \textbf{(P)+(M)}, where the right-hand-sides are of polynomial growth and fulfill \eqref{TriangStructure}. Furthermore, they have numerically investigated the global existence of solutions to the following system
			
			\begin{equation}\label{ExampleReactiondiffusionSystem2x2}\tag{$S_{\exp}$}
				\left\{
				\begin{array}{rcll}
					\partial_t u_1 + d_1(-\Delta)^{s}u_1 &=& -u_1 e^{(u_2)^\beta} & \text{in }(0,T)\times\Omega,\\
					\partial_t u_2 + d_2(-\Delta)^{s}u_2 &=& \phantom{-}u_1 e^{(u_2)^\beta} & \text{in }(0,T)\times\Omega,\\
					u_1=u_2 &=& 0 & \text{in }(0,T)\times(\mathbb{R}^N\setminus\Omega),\\
					u_1(0,\mathbf{x}) &=& u_{01}(\mathbf{x}), & \mathbf{x}\in\Omega,\\
					u_2(0,\mathbf{x}) &=& u_{02}(\mathbf{x}), & \mathbf{x}\in\Omega.
				\end{array}
				\right.
			\end{equation}

%				\begin{equation}\label{ExampleReactiondiffusionSystem2x2} \tag{$S_{\exp}$}
%					\left\{   \begin{array}{rclll}
%						\partial_t u_1(t,\textbf{x})+d_1(-\Delta)^{s}u_1(t,\textbf{x})&=&-u_1(t,\textbf{x})e^{(u_2(t,\textbf{x}))^\beta},&(t,\textbf{x})\in (0,T)\times\Omega,\\ 
%						\partial_t u_2(t,\textbf{x})+d_2(-\Delta)^{s}u_2(t,\textbf{x})&=&\phantom{-} u_1(t,\textbf{x})e^{(u_2(t,\textbf{x}))^\beta},&(t,\textbf{x})\in (0,T)\times\Omega,\\
%						u_1(t,\textbf{x})=u_2(t,\textbf{x})&=&0,&(t,\textbf{x})\in (0,T)\times(\R^N\setminus\Omega),\\
%						u_1(0,\textbf{x})&=&u_{01}(\textbf{x}),	& \textbf{x}\in\Omega,\\
%						u_2(0,\textbf{x})&=&u_{02}(\textbf{x}),	& \textbf{x}\in\Omega,
%					\end{array}    \right.
%				\end{equation}
			
			where $\beta\geq1$. Let us emphasize that this system arises in the modeling of exothermic  combustion in gases, see, $e.g.$, \cite{HerLacVel}.
			\\ (ii) In \cite{NguTang}, the authors addressed System \eqref{ReactiondiffusionSystem1} in the case  $\Omega=\R^N$.
		\end{itemize}
		\item[$\diamond$] \textbf{With spectral fractional Laplacians}
		\\
		In this case, as far as we know, the only available studies concern System \eqref{ReactiondiffusionSystem1} with $m=2$: \\
		(i) In \cite{SomBur2018}, the authors introduce a fractional model to describe coral growth in a tank, applied to modeling the pattern formation of coral reefs. More specifically, the model is characterized by the parameters $s_1= s_2$, $f_1(u_1,u_2)=1-u_1-\lambda u_1 u_2^2$ and $f_2(u_1,u_2)=-\mu u_2+\lambda u_1 u_2^2$, with $\lambda,\mu>0$.
		\\
		(ii) In \cite{ZhaoRuan}, the authors propose a fractional model to study  the long-range geographical spread of infectious diseases. To elaborate, the model is described by the parameters $s_1\neq s_2$, $f_1(u_1,u_2)=\varphi(\mathbf{x})u_1-\phi(\mathbf{x})u_1^2-\dfrac{\psi(\mathbf{x})u_1u_2}{u_1+u_2}+h(\mathbf{x})u_2$ and $f_2(u_1,u_2)=\dfrac{\psi(\mathbf{x})u_1u_2}{u_1+u_2}-h(\mathbf{x})u_2$, where the functions $\varphi$, $\phi$, $\psi$ and $h$ are H\"{o}lder continuous in $\overline{\Omega}$.
	\end{itemize}		
\end{itemize} 

To the best of our knowledge, System \eqref{ReactiondiffusionSystem1} with the reaction terms satisfying  {\bf (P)+(M)} and of polynomial growth has not yet been thoroughly investigated. In this work, we aim to address this by establishing the global existence of nonnegative strong solutions to System \eqref{ReactiondiffusionSystem1}. 
%This extends the study conducted in \cite{DaouLaamBaal2024}, where the author and her collaborators analyzed the same system with the regional fractional Laplacian ($RFL$) and $s_1=\cdots=s_m$.
This is the main purpose of this work. More precisely, we aim to generalize within the fractional setting ($0<s_i<1$) two main  well-established results  in the  classical case (\textit{i.e.} all $s_i$'s $=1$):
\begin{enumerate}
	\item 
	The first one concerns the typical case of reversible chemical reactions for three species, that is

%	\begin{equation}\label{SystemChemReac-DiffFR2-} 
%			\tag{$S_{\alpha,\beta,\gamma}$} 	
%			\left\{
%			\begin{array}{rclll}
%				\forall i=1,2,3,\hspace*{2.8cm}&&&\\
%				\partial_t u_i(t,\textbf{x})+d_i (-\Delta)_{Sp}^{s_i} u_i(t,\textbf{x})&=&	f_i(u_1(t,\textbf{x}),u_2(t,\textbf{x}),u_3(t,\textbf{x}) ), &(t,\textbf{x})\in (0,T)\times\Omega,\\ \mathcal{B}[u_i(t,\textbf{x})]&=&0,&(t,\textbf{x})\in (0,T)\times\partial\Omega, \\  u_i(0, \textbf{x})&=&u_{0i}(\textbf{x}),  & \textbf{x} \in \Omega,\end{array}\right.
%	\end{equation}

\begin{equation}\label{SystemChemReac-DiffFR2-}\tag{$S_{\alpha,\beta,\gamma}$}
	\left\{
	\begin{array}{rcll}
			\multicolumn{4}{l}{\forall i=1,2,3,} \\[2pt]
		\partial_t u_i + d_i(-\Delta)_{Sp}^{s_i}u_i
		&=& f_i(u_1,u_2,u_3), & \text{in }(0,T)\times\Omega,\\
		\mathcal{B}[u_i] &=& 0, & \text{in }(0,T)\times\partial\Omega,\\
		u_i(0,\mathbf{x}) &=& u_{0i}(\mathbf{x}), & \mathbf{x}\in\Omega.
	\end{array}
	\right.
\end{equation}

	 where $f_1= \alpha g,\; f_2= \beta g,\; f_3= -\gamma g \text{ with } g=u_3^{\gamma}-u_1^{\alpha}u_2^{\beta} \text{ and } 1\leq \alpha, \beta, \gamma <+\infty$.
	It is worth noting that this system naturally emerges in chemical kinetics when modeling the following reversible reaction
	\begin{equation}\label{reaction-chimique}
		\alpha U_1+\beta U_2 \rightleftharpoons \gamma U_3
	\end{equation}
	where $u_1,\, u_2,\,  u_3$ stand for the density of $U_1,\, U_2$ and $U_3$ respectively, {and $\alpha$, $\beta$, $\gamma$ are the stoichiometric coefficients}.
	\noindent  Moreover, let us point out that in  addition to satisfying {\bf (P)}, the functions $f_1$, $f_2$ and $f_3$ also fulfill {\bf (M)}, $i.e.$ $\beta\gamma f_1+\alpha\gamma f_2+2\alpha\beta f_3~=~0.$
	
	It should also be emphasized that global existence of strong nonnegative solutions to System \eqref{SystemChemReac-DiffFR2-}  in the classical case $s_1=s_2=s_3$ has been established in \cite{LaamActa2011}.

	\item  The second one addresses the case where System \eqref{ReactiondiffusionSystem1} consists of $m\geq 2$ equations and the $f_i$'s fulfill a {\it triangular structure} (see~{\eqref{TriangStructure}}). Let us mention that the classical version ($i.e.$ all $s_i$'s=1) of this problem has been studied in \cite{PierreSurvey2010}.
\end{enumerate}	

To establish our main theorems, we will adapt several well-known tools from the classical setting to the fractional framework. These include the maximal regularity theorem (see Theorem~\ref{BoundedSol}),  a Lamberton-type estimation in $L^p$ (see Theorem~\ref{DualProblemIneq}), the so-called Pierre's duality Lemma (see Lemma~\ref{FractionalDualityThm}). As emphasized in \cite{DaouLaamBaal2024}, these results are noteworthy  in themselves. 
That being said, it should be noted that  not all results from the classical case can be directly extended to the fractional setting; and when it is, it often entails substantial challenges.
To illustrate, in a different context, we refer to~{\cite{AtmaBirDaouLaam-Parabo,AtmaBirDaouLaam-Acta,AtmaBirDaouLaam-s1-s2-Ellip}}.

\medbreak

On the other hand, we will present some numerical simulations for System~\eqref{SystemChemReac-DiffFR2-}
in order to address an open theoretical question.
\\
\noindent
Let us recall that, in the classical diffusion case $s_i=1$, the global existence
and convergence towards an equilibrium are well understood
(see, e.g., \cite{PierSuzUma2018,FellLaam2016}).
More precisely, global strong solutions exist  when $\gamma>\alpha+\beta$ for bounded nonnegative initial data. However, the question of global existence is open whenever
$\gamma\leq\alpha+\beta$.  %under suitable hypothesis on initial data.
\\
In the case of spectral fractional Laplacians $0<s_i<1$, the global theory is
far from complete, with open questions arising both from the reaction parameters
$(\alpha,\beta,\gamma)$ and from the relative orders of diffusion. 
In particular, the global existence of strong solutions is still unknown when
$
s_3>\min\{s_1,s_2\}$ and $
\gamma\leq\alpha+\beta$.
In the numerical part of the paper, we focus precisely on this open configuration.

%\textcolor{olive}{Our purpose is to treat the case where at least two powers are different. The case $s_i=s_j$ is worth commenting. In this case, once we have the semigroup properties fulfilled, all theoretical results presented in our paper could be extended.}

\smallbreak
Before concluding the introduction, it is worth noting that the global existence of strong solutions to System $(S)$ with regional fractional Laplacians of different orders, will be addressed in a forthcoming work,  in the case where $\Omega$ is either bounded or equal to $\mathbb{R}^N$.

\smallbreak
The structure of this article is as follows.  In Section 2,
we provide the necessary background for the subsequent sections. 
First, we give two equivalent definitions of the $SFL$, the appropriate  framework and some fundamental properties of the corresponding semigroup. Then, we address a fractional evolution equation. More precisely, we investigate the criteria for ensuring the boundedness of the solution. Additionally, we focus on relevant concepts such as the dual problem, maximal regularity and the comparison principle.	
Section 3 is divided into four subsections. In subsection 3.1,   we begin by investigating the local existence of strong solutions to System \eqref{ReactiondiffusionSystem1}. Subsection 3.2 is dedicated to extending Pierre's duality lemma to the fractional case with different orders. 
Subsections 3.3 and 3.4 are devoted respectively to the proofs of our two main global existence theorems. In Section 4, we present some numerical simulations to explore the global existence of solutions to System \eqref{SystemChemReac-DiffFR2-} in an open theoretical case. 
	
	\medbreak	
	Now let us fix some notations that will be used throughout this work.
	\smallbreak 	\noindent{\bf Notations.}
	\begin{itemize}
		\item[---] $\Omega$ is a bounded regular open subset of $\mathbb{R}^N$ {with $N\in\mathbb{N}^*$}. 
		\item[---]  For any $T>0$,  $Q_T:= (0,T)\times \Omega$ and $\Sigma_T:= (0,T)\times \partial\Omega$.
		\item[---]For any $p\in[1,+\infty)$, $\|\phi\|_{L^p(\Omega)}=\left(\dis\int_\Omega|\phi(\mathbf{x})|^pd\mathbf{x} \right)^{\frac1p}$ and $\|\psi\|_{L^p(Q_T)}=\left(\dis\int_0^T\hskip-2mm\int_\Omega|\psi(t,\mathbf{x})|^p dt d\mathbf{x} \right)^{\frac1p}$.
		\item[---] $\|\phi\|_{L^\infty(\Omega)}= \esssup\limits_{\mathbf{x}\in\Omega} |\phi(\mathbf{x})|$ and  $\|\psi\|_{L^\infty(Q_T)}= \esssup\limits_{(t,\mathbf{x})\in Q_T} |\psi(t,\mathbf{x})|$.
		\item[---] The boundary operator $\mathcal{B}$ is defined as:
		$$
		\mathcal{B}[u]=u \quad\text{or}\quad \mathcal{B}[u]=\frac{\partial u}{\partial \nu} \quad \text {on} \quad \partial \Omega.
		$$		
\end{itemize}

\section{Preliminaries}
\label{sec:1}
In this section, we briefly review some results that will be relevant for later use. First, we recall the spectral fractional Laplacian and its semigroup.  Subsequently, we summarize existence and regularity results for solutions to a fractional evolution problem, along with a Lamberton-type estimate related to the dual problem. 
\subsection{Spectral fractional Laplacian ($SFL$)}
First, we set up the functional framework. Then, we recall two equivalent definitions of the spectral fractional Laplacian and discuss properties of the corresponding semigroup.
\subsubsection{Functional framework}
Let $\{\lambda_k\}_{k\in\mathbb{N}^*}$ be the  positive eigenvalues of $-\Delta$ in $\Omega$
with homogeneous Dirichlet or Neumann boundary conditions, ordered as
$0<\lambda_1\le \lambda_2\le\cdots\uparrow+\infty$.
Furthermore, let $\{e_k\}_{k\in \mathbb{N}^*}$ be the corresponding eigenfunctions.
Namely, for any $k\in \mathbb{N}^*$, $e_k$ is the solution to
$$
\left\{
\begin{array}{rclll}
	-\Delta e_k(\mathbf{x})&=&\lambda_k e_k(\mathbf{x}), & \mathbf{x}\in \Omega, \\ 
	\mathcal{B}[e_k(\mathbf{x})]&=&0, &  \mathbf{x}\in \partial\Omega.
\end{array}
\right.
$$
Also, the eigenfunctions are normalized in $L^2(\Omega)$, $i.e.$ $\dis\int_\Omega e_i(\mathbf{x}) e_j(\mathbf{x}) d\mathbf{x} =\delta_{ij}$, where $\delta_{ij}$ is the Kronecker symbol. In the Neumann case, there is an additional eigenvalue $\lambda_0=0$ associated with a constant eigenfunction normalized in $L^2(\Omega)$. For more details, see, for instance, \cite[Theorem 1.2.2 and Theorem~1.2.8]{Henrot2006}.

In this setting, we can write, for any $u\in L^2(\Omega)$, $u(\mathbf{x})=\displaystyle\sum\limits_{k=1}^{+\infty} u_k e_k(\mathbf{x})$ such that $u_k=\displaystyle\int_\Omega u(\mathbf{y}) e_k (\mathbf{y}) d\mathbf{y}$.

Moreover, for any $s>0$, let us introduce the following space
$$
\mathcal{X}^s(\Omega):= \{u\in L^{2}(\Omega)\;\;\text{such that}\;\; \sum\limits_{k=1}^{+\infty} \lambda_{k}^{s} u^{2}_k<{+\infty} \}.
$$
\noindent We easily verify that $(\mathcal{X}^s(\Omega), \langle\cdot,\cdot \rangle_{\mathcal{X}^s(\Omega)})$ is a Hilbert space, where
$$
\langle u,v\rangle_{\mathcal{X}^s(\Omega)}:=\sum\limits_{k=1}^{+\infty}  \lambda_{k}^{s} u_kv_k, \quad\text{for any }\; u,v\in \mathcal{X}^s(\Omega).
$$
For further details on this space, we refer the interested reader to, $e.g.$, \cite[Section 1]{HulVan1993} or \cite[Subsection 3.1]{BonSirVaz2015}.

\subsubsection{First definition of the $SFL$}

Let  $\mathbf{x}\in\Omega$ and $s\in(0,1)$. The spectral fractional Laplacian ({\em SFL}) is defined as
\begin{eqnarray}\label{DefinitionofSFL1}
	(-\Delta)_{Sp}^s u(\mathbf{x}):=\sum\limits_{k=1}^{+\infty} \lambda_{k}^su_{k}e_{k} (\mathbf{x}),
\end{eqnarray}
and its domain is given by
$$
\begin{array}{rcll}
	D(	(-\Delta)_{Sp}^s)&:=& \{u\in L^{2}(\Omega)\;\;\text{such that}\;\;  \|	(-\Delta)_{Sp}^su\|_{L^2(\Omega)}^2=\sum\limits_{k=1}^{+\infty} \lambda_{k}^{2s} u^{2}_k<{+\infty} \}
	\\
	&=& 
	%\mathcal{X}^{s}(\Omega)\cap\mathcal{X}^{2s}(\Omega)=
	\mathcal{X}^{2s}(\Omega).
\end{array}
$$
Obviously, when $s=0$ we recover the identity, and when 
$s=1$ we obtain the classical Laplacian.

\subsubsection{Semigroup properties} \label{SG}
Let us consider the operator $A^s$ on $L^2(\Omega)$, defined by
$$
\left\{
\begin{array}{rclll}
	D(A^s)&=&\mathcal{X}^{2s}(\Omega),
	\\
	A^su&=&-d(-\Delta)^s_{Sp}u,\;\;d>0,
\end{array}
\right.
$$
for any $u\in D(A^s)$.
It is well-known that $A^s$ generates a strongly continuous semigroup $\{T_s(t)\}_{t\geq 0}:=\{e^{tA^s}\}_{t\geq 0}$ on $L^2(\Omega)$. For readers who may be unfamiliar with the notion of semigroups, we suggest consulting \cite[Chapter 7]{Vrabie2003} or \cite[Chapter~1]{Pazy1983}. Furthermore, the semigroup $\{T_s(t)\}_{t\geq 0}$ is submarkovian (positivity-preserving and $L^\infty$-contractive), $namely$
$$T_s(t) u \geq 0\;\; \text{for any}\;\; t \geq 0, \;\;\text{as}\;\; 0 \leq u \in L^2(\Omega)$$ 
and
$$\|T_s(t) u\|_{L^p(\Omega)} \leq\|u\|_{L^p(\Omega)} \;\; \text{for any}\;\; u \in L^p(\Omega) \cap L^2(\Omega),\;\; p \in[1,+ \infty].$$
{As a consequence, $\{T_s(t)\}_{t\geq 0}$ can be extended to  a contraction semigroup on $L^p(\Omega)$ for any $p \in$ $[1,+\infty]$. The semigroup is strongly continuous if $p \in[1,+\infty)$ and bounded analytic if $p \in(1,+\infty)$ (see, for instance, \cite[Theorem 1.4.1]{Davies1990} for more details).} Throughout the paper, by a slight abuse of notation, we keep the same symbol
$T_s(t)$ for its extensions to $L^p(\Omega)$, and we denote again by $A^s$ its generator on $L^p(\Omega)$ when $p\in(1,+\infty)$.
In addition, $\{T_s(t)\}_{t\geq 0}$ is ultracontractive,  in the sense that  for any $1 \leq p \leq q \leq+ \infty$, the operator $T_s(t)$ maps $L^p(\Omega)$ into $L^q(\Omega)$. In other words,  there exists $C>0$ such that for any $u \in L^p(\Omega)$
\begin{equation}
	\label{Ultracontractive}
	\left\|T_s(t) u\right\|_{L^q(\Omega)} \leq C t^{-\frac{N}{2s}\left(\frac{1}{p}-\frac{1}{q}\right)}\|u\|_{L^p(\Omega)}, \quad \forall t>0 .
\end{equation}
For further details and proofs regarding the properties of $\{T_s(t)\}_{t\geq 0}$ , see \cite[Chapter 2]{GalWarma2020}.
\subsubsection{Second definition of the $SFL$}

Let $s\in(0,1)$. The operator {\em SFL} can also be expressed as an integral with respect to an appropriate kernel, along with an additional linear term, see \cite[Definition 1 and Lemma 38]{AbaDup2017} and \cite[Subsection 2.3]{AbaVal2019}.  More precisely,  for any $u\in\mathcal{X}^{2s}(\Omega)$ and for $a.e.$ $\mathbf{x}\in\Omega$, we have
\begin{equation}\label{DefinitionofSFL2}
	(-\Delta)_{Sp}^s u(\mathbf{x})=\text{P.V.} \int_{\Omega}[u(\mathbf{x})-u(\mathbf{y})] J(\mathbf{x}, \mathbf{y}) d \mathbf{y}+\kappa(\mathbf{x}) u(\mathbf{x}),
\end{equation}
where P.V. stands for the Cauchy principal value,
$$ J(\mathbf{x}, \mathbf{y}):=\frac{s}{\Gamma(1-s)} \int_0^{+\infty} \frac{p_{\Omega}(t, \mathbf{x}, \mathbf{y})}{t^{1+s}} d t \geq 0,$$
and $p_{\Omega}(t, \mathbf{x}, \mathbf{y}) $ is the heat kernel of the Dirichlet or Neumann Laplacian. Furthermore,
$$
\kappa(\mathbf{x}):=\frac{s}{\Gamma(1-s)} \int_0^{+\infty}\left(1-\int_{\Omega} p_{\Omega}(t, \mathbf{x}, \mathbf{y}) d \mathbf{y}\right) \frac{d t}{t^{1+s}}\geq 0,
$$
in the case of the Dirichlet Laplacian, and $\kappa (\mathbf{x})=0$ in the case of the Neumann one. For more details, see \cite[Subsection 1.1 and Subsection 1.2]{DipGiaVal2022}.

\subsection{Evolution Problem governed by the $SFL$ and its dual}
In this subsection, we review some existence results and key properties of solutions to the fractional evolution problem \eqref{ParabolicProblem}. Furthermore, we examine the dual of this problem to derive a useful Lamberton-type estimate.

\subsubsection{Fractional evolution problem}
First, let us consider the following problem:
\begin{equation}\label{ParabolicProblem}
	\left\{\begin{array}{rclll}
		\partial_t w(t,\textbf{x})+d (-\Delta)_{Sp}^s w(t,\textbf{x})&=&h(t,\textbf{x}), & (t,\textbf{x})\in Q_T, \\
		\mathcal{B}[w(t,\textbf{x})]&=&0, & (t,\mathbf{x})\in\Sigma_T, \\
		w( 0,\textbf{x})&=&w_{0}(\textbf{x})\geq 0, & \textbf{x}\in\Omega,
	\end{array}\right.
\end{equation}
where $s\in(0,1)$, $d>0$, $w_0\in L^p(\Omega)$ and $h\in L^p(Q_T$) with $p\geq1$. 
\vskip2mm 
\begin{definition}[Weak solution] \label{Weak_solution}
	Let $h\in L^1(0,T; L^p(\Omega))$ with $p>1$. A function $ w\in \mathcal{C}([0,T];L^p(\Omega))$ is called a \textit{weak} solution to Problem \eqref{ParabolicProblem} if for any $(t,\mathbf{x})\in[0,T)\times\Omega$,  
	\begin{equation}\label{weaksolution}
		w(t,\mathbf{x})=T_s(t) w_0(\mathbf{x})+\int_0^t T_s(t-\tau) h(\tau,\mathbf{x}) d \tau,
	\end{equation}
	where $\{T_s(t)\}_{t\geq 0}$ is the semigroup generated by $Azx^s$ for any $s\in(0,1)$ (see Paragraph~\ref{SG}).
\end{definition}

\begin{definition}[Strong solution]\label{Strong_solution}
	Let $p>1$ and $s\in (0,1)$.	A function $ w$ is said to be a strong solution to Problem  \eqref{ParabolicProblem}  %on $[0, T)$ 
	if :
	
	{\em (i)} $ w \in \mathcal{C}([0,T);L^p(\Omega))\cap\mathcal{C}^1(0,T;L^p(\Omega))$ ;
	
	{\em (ii)} for any $t\in(0,T)$, $ w(t,\cdot)\in D(A^s) $ ; 
	
	{\em (iii)} Problem \eqref{ParabolicProblem} is satisfied $a.e.$ in $Q_T$.
\end{definition}

The following theorem guarantees the existence of  a unique weak solution.
\begin{theorem}[{\cite[Chapter 4]{Pazy1983}}]\label{PazyExistence}
	Let $h\in L^1(0,T; L^p(\Omega))$ and $w_0\in L^p(\Omega)$ with $p>1$. Then, Problem \eqref{ParabolicProblem} admits a unique weak solution.
\end{theorem}

\begin{remark}
	Any ``weak solution'' belonging to $L^\infty(Q_T)$ is, in fact,
	regular enough to be a ``strong solution''.
\end{remark}

In the three subsequent theorems, we give some regularity properties of the weak solution to Problem \eqref{ParabolicProblem}.

\medbreak
\noindent {1)} A maximal regularity result in $L^\infty$:

\begin{theorem} \label{BoundedSol}%Theorem 3.1
	Let $s\in(0,1)$, $w_0 \in L^\infty(\Omega)^+$ and $h\in L^p(Q_T$) with $p>1$. Moreover, let $w$ be the weak solution to Problem \eqref{ParabolicProblem}. Then, for any $p>\dfrac{N+2s}{2s}$ we have
	$$
	\|w\|_{L^\infty(Q_T)}\leq \|w_0\|_{L^\infty(\Omega)}+C\|h\|_{L^p(Q_T)}, \quad \text{for some } \; C>0.
	$$
\end{theorem}
\begin{proof}
	By applying the ultracontractivity estimate \eqref{Ultracontractive}, the proof follows as in \cite[Theorem 3.1]{DaouLaamBaal2024}.
\end{proof}

{
	\noindent{2)} A first comparison principle result :
}

\begin{theorem}\label{maximumprinciple1} 	Let $w_0 \in L^\infty(\Omega)^+$  and $h\in L^p(Q_T)$ with $p>1$. Further, let $w$ be the weak solution to Problem \eqref{ParabolicProblem}. If moreover $h\geq 0$ $a.e.$ in $Q_T$, then $w\geq 0$ $a.e.$ in $Q_T$. 
\end{theorem}
\begin{proof} Using the definition \eqref{DefinitionofSFL2}, the proof proceeds similarly to \cite[Theorem 3.2]{DaouLaamBaal2024}.
\end{proof}
{
	\noindent{3)} A second comparison principle result :
}
\begin{theorem}\label{maximumprinciple2} 	Let $w_0 \in L^\infty(\Omega)^+$, $h\in L^p(Q_T$) with $p>1$, and $w$ be the weak solution to Problem \eqref{ParabolicProblem}. In addition, let $z$ be the weak solution to
	\begin{equation}\label{ParabolicProblem2}
		\left\{\begin{array}{rclll}
			\partial_t z(t,\mathbf{x})+d (-\Delta)_{Sp}^s z(t,\mathbf{x})&=&g(t,\mathbf{x}), & (t,\mathbf{x})\in Q_T, \\
			\mathcal{B}[z(t,\mathbf{x})]&=&0, & (t,\mathbf{x})\in\Sigma_T, \\
			z( 0,\mathbf{x})&=&z_{0}(\mathbf{x})\geq 0, & \mathbf{x}\in\Omega,
		\end{array}\right.
	\end{equation}
	where  $s\in(0,1)$, $z_0 \in L^\infty(\Omega)$ and	$g\in L^q(Q_T$) for some $q>1$. Suppose that $h \leq g$ and $w_{0} \leq z_{0}$ $a.e.$ in $Q_T$ and $\Omega$ respectively. Then $w(t,\mathbf{x}) \leq z(t,\mathbf{x})$ for $a.e.$ in $(t,\mathbf{x})\in Q_T$.
\end{theorem}
\begin{proof}
	Let us consider the following problem
	\begin{equation}\label{ParabolicProblem3}
		\left\{\begin{array}{rclll}
			\partial_t (z-w)(t,\mathbf{x})+d (-\Delta)_{Sp}^s (z-w)(t,\mathbf{x})&=&(g-h)(t,\mathbf{x}), & (t,\mathbf{x})\in Q_T, \\
			\mathcal{B}[(z-w)(t,\mathbf{x})]&=&0, & (t,\mathbf{x})\in\Sigma_T, \\
			(z-w)( 0,\mathbf{x})&=&(z_{0}-w_{0})(\mathbf{x}), & \mathbf{x}\in\Omega.
		\end{array}\right.
	\end{equation}
	To achieve the proof, we apply Theorem \ref{maximumprinciple1} to $z-w$.
\end{proof}

	\subsubsection{Dual of Problem \eqref{ParabolicProblem}}

The dual of Problem \eqref{ParabolicProblem} is given by
\begin{equation}\tag{$P_{\phi,T}$}\label{DualProblem}
	\left\{\begin{array}{rclll}
		-\partial_t \mathcal{V}(t,\textbf{x})+d(-\Delta)_{Sp}^s \mathcal{V}(t,\textbf{x})&=&\phi(t,\textbf{x}), & (t,\textbf{x})\in Q_T, \\
		\mathcal{B}[\mathcal{V}(t,\textbf{x})]&=&0, & (t,\textbf{x})\in \Sigma_T, \\
		\mathcal{V}(T, \textbf{x})&=& 0, & \textbf{x}\in\Omega,
	\end{array}\right.
\end{equation}
where $s\in(0,1)$, $d>0$ and
$\phi$  is a regular function.
\medbreak

\noindent The following theorem outlines the conditions for the existence of a solution to Problem \eqref{DualProblem},  and ensures it meets a  Lamberton-type estimation in $L^p$.

\begin{theorem} \label{DualProblemIneq} %Theorem 3.2
	Let  $s\in(0,1)$.	Assume that $\phi\in L^p(Q_T)$ with $p>1$. Then, Problem 
	\eqref{DualProblem} admits a unique weak solution. Moreover,
	there exists a constant $C:=C(p,T)>0$ such that    	\begin{equation}\label{DualityInequality}
		\left\| \partial_t \mathcal{V}\right\|_{L^p(Q_T)}+  \|\mathcal{V}\|_{L^p(Q_T)}+d \left\|(-\Delta)_{Sp}^s\mathcal{V}\right\|_{L^p(Q_T)}+\|\mathcal{V}_0\|_{L^p(\Omega)}\leq C   \left\|\phi\right\|_{L^p(Q_T)},
	\end{equation}
	where $\mathcal{V}_0:=\mathcal{V}(0,\cdot)$.
\end{theorem}   
\begin{proof}
	By leveraging the properties of the semigroup $\{T_s(t)\}_{t \geq 0}$, the proof can be constructed by following the steps of the proof of \cite[Theorem 3.3]{DaouLaamBaal2024},  relying on \cite[Theorem 1]{Lamberton1987}.
\end{proof}

\section{Global existence for reaction-diffusion systems governed by spectral fractional Laplacians}
Let us recall our main system:
\begin{equation}\label{ReactiondiffusionSystem3}\tag{$S$}
	\left\{   \begin{array}{rcll}
			\multicolumn{4}{l}{ \forall i=1,\ldots,m,} \\[2pt]
		\dis \partial_t u_i(t,\textbf{x})+d_i(-\Delta)_{Sp}^{s_i}u_i(t,\textbf{x})&=&f_i(u_1(t,\textbf{x}),\ldots,u_m(t,\textbf{x})),&(t,\textbf{x})\in Q_T,\\ \dis
		\mathcal{B}[u_i(t, \textbf{x})]&=&0,&(t,\textbf{x})\in \Sigma_T,\\\dis
		u_i(0,\textbf{x})&=&u_{0i}(\textbf{x}),	&\textbf{x}\in\Omega,
	\end{array}
	\right.
\end{equation}
where  $m\geq 2$ and for each $i\in\llbracket1, m\rrbracket$, $0<s_i<1$, $d_i>0$ and $f_i$ is locally Lipschitz continuous. 

\smallbreak
As stated in the introduction, our main contribution is to extend two fundamental results from the classical case to the fractional setting. This is the aim of this section. For clarity, it is structured into four subsections.  In the first one, we establish the local existence of strong solutions to System \eqref{ReactiondiffusionSystem3}. In the second, we prove a result that enables the extension of Pierre's duality lemma to the fractional case with different orders, and we subsequently demonstrate this extension. Finally, Subsections 4.2 and 4.3 are respectively devoted to the proofs of our main global existence theorems.

\subsection{\bf Local existence}
{The local existence lemma is well-established. Here, we will recall its statement and omit the proof, as it is straightforward.
	\begin{lemma}\label{LocalExistence}
		Let
		$(u_{01},\ldots,u_{0m}) \in (L^\infty(\Omega))^m$. Assume that the $f_i$'s  are  locally Lipschitz continuous. 
		
		Then,
		there exists $T_{\max}>0$ and 
		$(\varphi_1,\ldots,\varphi_m)\in \mathcal{C}([0,T_{\max}),[0, +\infty)^m)$ such that :\\
		{\em (i)} System~\eqref{ReactiondiffusionSystem1} admits  a unique strong solution $(u_{1},\ldots,u_{m})$ in $Q_{T_{\max}}$; \\
		{\em (ii)}	for each $i\in \llbracket1,m\rrbracket$
		\begin{equation}
			\| u_i(t,\cdot)\|_{L^\infty(\Omega)}\leq \varphi_i(t)\;\; \text{for any}\;\; t\in(0,T_{\max});
		\end{equation}
		{\em (iii)}	if $T_{\max}<+\infty$, then $\dis\lim_{t\nearrow T_{\max}}\dis\sum\limits_{i=1}^m \|u_i(t,\cdot)\|_{L^\infty(\Omega)}=+\infty$;\\
		{\em (iv)}	if, in addition,  the $f_i$'s satisfy {\bf (P)}, then 
		$$\big(\forall i \in\llbracket1,m\rrbracket,\; u_{0i}(.)\geq 0\big) 	\Longrightarrow \big(\forall i \in\llbracket1,m\rrbracket,\; u_i(t,.)\geq 0 \; \forall t\in [0,T_{\max})\big).$$
	\end{lemma}
	\begin{remark}
		In light of {\em (iii)}, to prove global existence of strong solutions to System~\eqref{ReactiondiffusionSystem1},  it suffices to establish  an a priori estimate of the form:
		\begin{equation}\label{crucial-estimate}
			\forall t\in [0,T_{\max}),\; \dis\sum\limits_{i=1}^m \|u_i(t,\cdot) \|_{L^\infty(\Omega)} \leq\psi(t),  
		\end{equation}
		where $\psi :[0,+\infty)\rightarrow [0,+\infty)$ is a continuous function.
	\end{remark}
	\begin{remark}
		It turns out that establishing an estimate like \eqref{crucial-estimate} is far from direct, except in the trivial case \( d_1 = \cdots = d_m \), \( s_1 = \cdots = s_m \), and the \( f_i \)'s fulfill {({\bf P})+({\bf M})}. For this particular instance, System \eqref{ReactiondiffusionSystem3} admits a unique nonnegative global strong solution regardless of the growth of the \( f_i \)’s. The proof can be directly constructed following \cite[Subsection 4.1]{DaouLaamBaal2024}.  Nonetheless, the situation becomes considerably more complicated when the diffusion coefficients differ. Additional assumptions on the growth of the source terms are required to guarantee the global existence of strong solutions, even in the case \( s_1 = \cdots = s_m = 1 \).  Moreover, the complexity increases further when different fractional orders are taken into account.
	\end{remark}
	\subsection{\bf Fractional version of Pierre's duality lemma}
	To prove our global existence theorems, we will rely heavily on the following lemma. Let us mention that when the fractional powers tend to 1, the following result  is essentially a generalization of Pierre's duality lemma (see \cite{HolMarPie1987,PierreSurvey2010}  and references therein). 
	\begin{lemma}\label{FractionalDualityThm}
		
		Let $T>0$ and $s_1,s_2 \in (0,1)$ with $s_1\leq s_2$.
		Consider the following $2\times2$ system :
		\begin{equation}\label{ReactiondiffusionSystem2}\tag{$S_{2\times2}$}
			\left\{   \begin{array}{rcll}
				\dis \partial_t u_1(t,\mathbf{x})+d_1(-\Delta)_{Sp}^{s_1}u_1(t,\mathbf{x})&=&f_1(u_1(t,\mathbf{x}),u_2(t,\mathbf{x})),&(t,\mathbf{x})\in Q_T,
				\\
				\dis \partial_t u_2(t,\mathbf{x})+d_2(-\Delta)_{Sp}^{s_2}u_2(t,\mathbf{x})&=&f_2(u_1(t,\mathbf{x}),u_2(t,\mathbf{x})),&(t,\mathbf{x})\in Q_T,
				\\ \dis
				\mathcal{B}[u_1(t, \mathbf{x})]=\mathcal{B}[u_2(t, \mathbf{x})]&=&0,&(t,\mathbf{x})\in \Sigma_T,\\\dis
				u_1(0,\mathbf{x})&=&u_{01}(\mathbf{x}),	&\mathbf{x}\in\Omega,\\\dis
				u_2(0,\mathbf{x})&=&u_{02}(\mathbf{x}),	&\mathbf{x}\in\Omega,
			\end{array}
			\right.
		\end{equation}
		where $f_1$ and $f_2$ satisfy $(\mathbf{P})${\em +}$(\mathbf{M})$. Let $(u_1,u_2)$ be the strong solution to System~\eqref{ReactiondiffusionSystem2} in $Q_T$. Then, for any $p\in(1,+\infty)$, there exists  $\lambda > 0$ such that
		\begin{equation}\label{PierreEstimate}
			\|u_j\|_{L^p\left(Q_T\right)} \leq \lambda \left(1+\|u_i\|_{L^p\left(Q_T\right)}\right),
		\end{equation}
		where $ i,j\in \{1,2\}$ with $i\neq j$ and $s_i\leq s_j$.
	\end{lemma}
	
	\medbreak
	
	{Before giving the proof of Lemma \ref{FractionalDualityThm}, let us first state and prove a useful estimate that will be employed in the subsequent.

	\begin{proposition}\label{estimation}
		Let $p\in(1,+\infty)$ and $s_1,s_2\in(0,1)$ with $s_1\le s_2$. Then there exists
		$C>0$ such that for any
	$u\in D(A^{s_2})\subset L^p(\Omega)$,
		\begin{equation}\label{inequality_s1<s2_correct}
			\|(-\Delta)_{Sp}^{s_1}u\|_{L^p(\Omega)}
			\le C\,\|(-\Delta)_{Sp}^{s_2}u\|_{L^p(\Omega)}^{\frac{s_1}{s_2}}
			\|u\|_{L^p(\Omega)}^{\frac{s_2-s_1}{s_2}}.
		\end{equation}
	\end{proposition}
	
	\begin{proof}
		To start, the case $s_1=s_2$ is trivial. Now, let us assume that $s_1<s_2$.
		Set $\theta:=\frac{s_1}{s_2}\in(0,1)$ and define the operator
		\[
		A_0:=(-\Delta)_{Sp}^{s_2}\quad\text{on }L^p(\Omega),
		\]
		where $A_0$ is understood as the generator on $L^p(\Omega)$ associated with the
		$L^p$-extension of the semigroup (see Subsection~\ref{SG}). Since $-A_0$ generates
		a bounded analytic semigroup on $L^p(\Omega)$, the operator $A_0$ is 
		sectorial on $L^p(\Omega)$. Hence we may apply the moment inequality for
		fractional powers (see, e.g., \cite[Corollary~7.2]{Haase2005})  to obtain
		\begin{equation}\label{haase_moment}
			\|A_0^\theta u\|_{L^p(\Omega)}
			\le C\,\|A_0u\|_{L^p(\Omega)}^{\theta}\,\|u\|_{L^p(\Omega)}^{1-\theta},
			\qquad \forall\,u\in D(A_0),
		\end{equation}
		for some constant $C_0>0$.
	
	Finally, by the multiplicative property of fractional powers for sectorial
		operators (see, e.g., \cite[Formula~2.6]{Komatsu1967}),
		\[
		A_0^\theta=\bigl((-\Delta)_{Sp}^{s_2}\bigr)^{s_1/s_2}=(-\Delta)_{Sp}^{s_1},
		\]
		and $A_0u=(-\Delta)_{Sp}^{s_2}u$. Plugging these identities into
		\eqref{haase_moment} and using $1-\theta=\dfrac{s_2-s_1}{s_2}$ yields
		\eqref{inequality_s1<s2_correct}.
	\end{proof}

\begin{proof}[Proof of Lemma \ref{FractionalDualityThm}]  Let $T\in(0,T_{\max})$ and $(t,\mathbf{x})\in (0,T]\times\Omega$. Furthermore, let  $ i,j\in \{1,2\}$ with $i\neq j$ and assume that $0<s_i\leq s_j<1$.
	By hypothesis,  $f_i$ and $f_j$ satisfy {\bf (M)}, then there exists $a_i,a_j\in (0,+\infty)$ such that $a_i f_i(t,\mathbf{x})+a_j f_j(t,\mathbf{x})\leq 0$.
	By multiplying the $i$-th (resp. the $j$-th) equation of System \eqref{ReactiondiffusionSystem2} by $a_i$ (resp. $a_j$) and summing the two equations, we get for any $(t,\mathbf{x})\in Q_T$
	\begin{equation}\label{inequ_1+u_2}
		\partial_t(a_iu_i(t,\mathbf{x})+a_ju_j(t,\mathbf{x})) +a_id_i (-\Delta)_{Sp}^{s_i}u_i(t,\mathbf{x})+a_jd_j(-\Delta)_{Sp}^{s_j}u_j(t,\mathbf{x})\leq 0.
	\end{equation}
	Now, let $\phi$ be a nonnegative regular function and let $\mathcal{V}$ be the solution to Problem \eqref{DualProblem} with $s=s_j$ and $d=d_j$. First, let us multiply \eqref{inequ_1+u_2} by $\mathcal{V}$ and integrate over $Q_T$. Then, we obtain 
	\begin{equation}
		\left\{
		\begin{array}{ll}
			\dis
			-\iint_{Q_T} (a_iu_i(t,\mathbf{x})+a_ju_j(t,\mathbf{x})) \partial_t \mathcal{V}(t,\mathbf{x}) d\mathbf{x}dt\\\\\dis\hskip2cm+a_jd_j \iint_{Q_T} u_j(t,\mathbf{x}) (-\Delta)_{Sp}^{s_j} \mathcal{V}(t,\mathbf{x}) d\mathbf{x}dt \\\\\dis \leq
			\int_\Omega (a_iu_{0i}+a_ju_{0j})(\mathbf{x}) \mathcal{V}_0(\mathbf{x})d\mathbf{x}-a_id_i \iint_{Q_T} u_i(t,\mathbf{x}) (-\Delta)_{Sp}^{s_i}\mathcal{V}(t,\mathbf{x}) d\mathbf{x}dt .
		\end{array}
		\right.
	\end{equation}	 
	Then, by combining the terms, we get
	\begin{equation}
		\left\{
		\begin{array}{ll}
			\dis
			a_j\iint_{Q_T} u_j(t,\mathbf{x}) \big[ -\partial_t \mathcal{V}(t,\mathbf{x})+d_j (-\Delta)_{Sp}^{s_j}\mathcal{V}(t,\mathbf{x})  \big] d\mathbf{x}dt \\\\\dis\hskip2cm- a_i \iint_{Q_T}
			u_i(t,\mathbf{x}) \partial_t \mathcal{V}(t,\mathbf{x})  d\mathbf{x}dt \\\\\dis \leq
			\int_\Omega (a_iu_{0i}+a_ju_{0j})(\mathbf{x}) \mathcal{V}_0(\mathbf{x})d\mathbf{x}-a_id_i \iint_{Q_T} u_i(t,\mathbf{x}) (-\Delta)_{Sp}^{s_i}\mathcal{V}(t,\mathbf{x}) d\mathbf{x}dt .
		\end{array}
		\right.
	\end{equation}	 
	Hence,
	\begin{equation} \label{EstimateLemPierre}
		\left\{
		\begin{array}{ll}
			\dis a_j	\iint_{Q_T} u_j(t,\mathbf{x}) \phi(t,\mathbf{x}) d\mathbf{x}dt
			\\\\\dis
			\leq \int_\Omega (a_iu_{0i}+a_ju_{0j})(\mathbf{x}) \mathcal{V}_0(\mathbf{x})d\mathbf{x} \\\\\dis\hskip2cm+ a_i \iint_{Q_T}
			u_i(t,\mathbf{x})\big[ \partial_t \mathcal{V}(t,\mathbf{x})-d_i(-\Delta)_{Sp}^{s_i} \mathcal{V}(t,\mathbf{x})  \big] d\mathbf{x}dt .
		\end{array}
		\right.
	\end{equation}
	Now, let $p>1$. By applying H\"{o}lder's inequality  and using \eqref{DualityInequality},  we have
	\begin{equation}
		\int_\Omega (a_iu_{0i}+a_ju_{0j})(\mathbf{x}) \mathcal{V}_0(\mathbf{x})d\mathbf{x} \leq a_0C \|u_{0i}+u_{0j}\|_{L^p(\Omega)} \|\phi\|_{L^{p^\prime}(Q_T)},
	\end{equation}
	where $a_0:=\max\{a_i,a_j\}$ and $p^\prime:=\dfrac{p}{p-1}$. Moreover,
	\begin{equation}
		\left\{
		\begin{array}{ll}\dis
		a_i\iint_{Q_T}
		u_i(t,\mathbf{x})\big[ \partial_t \mathcal{V}(t,\mathbf{x})-d_i(-\Delta)_{Sp}^{s_i}\mathcal{V}(t,\mathbf{x})  \big] d\mathbf{x}dt \\\\\dis\leq a_i 
		\|u_i\|_{L^p(Q_T)} \left(C	\|\phi\|_{L^{p^\prime}(Q_T)}+d_i\|(-\Delta)_{Sp}^{s_i}\mathcal{V}\|_{L^{p^\prime}(Q_T)}\right).
			\end{array}
		\right.
	\end{equation}
	\\
	As $s_i\leq s_j$, using Proposition \ref{estimation},  we get
	$$
	\left\{
	\begin{array}{ll}\dis 
	\|(-\Delta)_{Sp}^{s_i}\mathcal{V}\|_{L^{p^\prime}(Q_T)}&\dis\leq C \|(-\Delta)_{Sp}^{s_j}\mathcal{V}\|_{L^{p^\prime}(Q_T)}^{\frac{s_i}{s_j}}  \|\mathcal{V}\|_{L^{p^\prime}(Q_T)}^{\frac{s_j-s_i}{s_j}}
	\\\\&\dis \leq C \|\phi\|^{\frac{s_i}{s_j}} _{L^{p^\prime}(Q_T)} \|\phi\|^{\frac{s_j-s_i}{s_j}} _{L^{p^\prime}(Q_T)}=C \|\phi\| _{L^{p^\prime}(Q_T)},
		\end{array}
	\right.
	$$
	where the last inequality follows by using again \eqref{DualityInequality}.
	Going back to \eqref{EstimateLemPierre}, we obtain
	\begin{equation}
			\left\{
		\begin{array}{ll}\dis
		a_j\iint_{Q_T} u_j(t,\mathbf{x}) \phi(t,\mathbf{x}) d\mathbf{x}dt \\\\\dis\leq C\big[ a_0 \|u_{0i}+u_{0j}\|_{L^p(\Omega)} \|\phi\|_{L^{p^\prime}(Q_T)} + a_i (1+d_i) 
		\|u_i\|_{L^p(Q_T)} 	\|\phi\|_{L^{p^\prime}(Q_T)}\big].
			\end{array}
		\right.
	\end{equation}
	\noindent Lastly, we deduce by duality that
	$$
	\|u_j\|_{L^p(Q_T)} \leq \frac{C}{a_j}\left(a_0\left\|u_{0i}+u_{0j}\right\|_{L^p(\Omega)}+a_i\left(1+d_i\right)\|u_i\|_{L^p(Q_T)}\right),
	$$
	which implies
	$$
	\|u_j\|_{L^p(Q_T)} \leq \frac{C}{a_j} \max \left\{a_0\left \|u_{0i}+u_{0j}\right\|_{L^p(\Omega)},a_i(1+d_i)\right\}\left(1+\|u_i\|_{L^p\left(Q_T\right)}\right).
	$$
	Accordingly, we get the desired estimation.
\end{proof}

	\subsection{Global existence : case of reversible chemical reaction with three species}\label{reversible-chemical-reaction}
In this section,  we will state and prove our first global existence related to the following system:
\begin{equation}\label{SystemChemReac-DiffFR2} 
	\tag{$S_{\alpha,\beta,\gamma}$} 	
	\left\{
	\begin{array}{rclll}
		\multicolumn{4}{l}{\forall i=1,2,3,} \\[2pt]
		\partial_t u_i(t,\textbf{x})+d_i (-\Delta)_{Sp}^{s_i} u_i(t,\textbf{x})&=&	f_i(u_1(t,\textbf{x}),u_2(t,\textbf{x}),u_3(t,\textbf{x}) ), &(t,\textbf{x})\in Q_T,\\ \mathcal{B}[u_i(t,\textbf{x})]&=&0,&(t,\textbf{x})\in \Sigma_T, \\  u_i(0, \textbf{x})&=&u_{0i}(\textbf{x}),  & \textbf{x} \in \Omega,\end{array}\right.
\end{equation}

\noindent where $0<s_i<1$, $f_1= \alpha g,\; f_2= \beta g,\; f_3= -\gamma g \text{ with } g=u_3^{\gamma}-u_1^{\alpha}u_2^{\beta} \text{ and } 1\leq \alpha, \beta, \gamma <+\infty$. 	

\noindent
Let us mention that for $\left(r_1, r_2, r_3\right) \in\left[0,+\infty\left){ }^3\right.\right.$, we have

$$
f_1\left(0, r_2, r_3\right)=\alpha r_3^\gamma \geq 0, \; f_2\left(r_1, 0, r_3\right)=\beta r_3^\gamma \geq 0 \; \text { and } f_3\left(r_1, r_2, 0\right)=\gamma r_1^\alpha r_2^\beta \geq 0
$$
and

$$
\beta \gamma f_1\left(r_1, r_2, r_3\right)+\alpha \gamma f_2\left(r_1, r_2, r_3\right)+2 \alpha \beta f_3\left(r_1, r_2, r_3\right)=0.
$$
Then, System \eqref{SystemChemReac-DiffFR2} fulfills $(\mathbf{P})$ and $(\mathbf{M})$.

\smallbreak
Our first result of global existence is as follows:

\begin{theorem}  \label{MainTh1}%\label{Thm-ReversibleChemicalReaction-DiffCL}%
	For each $i\in\llbracket1,3\rrbracket$, let $s_i\in(0,1)$ and
	$u_{0i} \in L^\infty(\Omega)^+$. Moreover, let $(\alpha, \beta, \gamma)\in  [1,+\infty)^3$. Then, System~\eqref{SystemChemReac-DiffFR2} admits a unique nonnegative global strong solution if one of the following assumptions holds: 
	\\
	{\em (i)}  $s_3\leq\min\{s_1,s_2\}$ and $\gamma>\alpha+\beta$ ;\\
	{\em (ii)} $s_3\geq\max\{s_1,s_2\} $ and $(\alpha, \beta,\gamma)\in  [1,+\infty)^2\times \{1\}$.
\end{theorem}

{\bf Comment.} It is noteworthy that the global existence of strong solutions when $2<\gamma\leq\alpha+\beta$ has remained an open question since 2011, even in the case of the classical Laplacian, as discussed in \cite{LaamActa2011}.
This is particularly surprising from a chemical perspective, as the reaction is expected to be reversible, suggesting that a strong solution should exist globally, similar to the case $\alpha+\beta< \gamma$.

\begin{proof}
	By Lemma \ref{LocalExistence}, System \eqref{SystemChemReac-DiffFR2} admits a unique nonnegative strong solution $(u_1,u_2,u_3)$ in $Q_{T_{\max}}$.  It remains to prove that this solution is global, $i.e.$ $T_{\max}=+\infty$. We closely follow the proofs of \cite[Theorems 1 and 3]{LaamActa2011} and \cite[Theorem 4.1]{DaouLaamBaal2024}. Now, for any $T\in(0,T_{\max})$ and  $t\in (0,T]$, we proceed to prove  (i) and (ii).
	\medbreak
	\noindent{\bf(i)} Let $s_3\leq \min\{s_1,s_2\}$ and $(\alpha, \beta, \gamma)\in  [1,+\infty)^3$ such that $\gamma>\alpha+\beta$. Moreover, let us consider the two following problems
	
	\begin{equation}
		\label{P1} \tag{$P'_{\gamma,1}$} 
		\left\{	\begin{array}{rclll}\partial_t v (t,\mathbf{x})+d_1 (-\Delta)_{Sp}^{s_1} v(t,\mathbf{x})&=&\alpha(u_3(t,\mathbf{x}))^\gamma, & (t, \mathbf{x}) \in Q_T, \\ \mathcal{B}[v(t, \mathbf{x})]&=&0, & (t, \mathbf{x}) \in \Sigma_T, \\ v(0, \mathbf{x})&=&u_{01}(\mathbf{x}), & \mathbf{x} \in \Omega,\end{array}\right.
	\end{equation}
	and
	\begin{equation}
		\label{P2} \tag{$P'_{\gamma,2}$} 
		\left\{	\begin{array}{rclll}\partial_t z (t,\mathbf{x})+d_2 (-\Delta)_{Sp}^{s_2} z(t,\mathbf{x})&=&\beta(u_3(t,\mathbf{x}))^\gamma, & (t, \mathbf{x}) \in Q_T, \\ \mathcal{B}[z(t, \mathbf{x})]&=&0, & (t, \mathbf{x}) \in \Sigma_T, \\ z(0, \mathbf{x})&=&u_{02}(\mathbf{x}), & \mathbf{x} \in \Omega.\end{array}\right.
	\end{equation}
	Since $u_1$, $u_2$ and $u_3$ are nonnegative,  the right-hand sides of the equations for $i=1,2$ of System \eqref{SystemChemReac-DiffFR2} fulfill $$f_1(u_1(t,\textbf{x}),u_2(t,\textbf{x}),u_3(t,\textbf{x}))\leq\alpha (u_3(t,\mathbf{x}))^\gamma$$ and $$f_2(u_1(t,\textbf{x}),u_2(t,\textbf{x}),u_3(t,\textbf{x}))\leq \beta(u_3(t,\mathbf{x}))^\gamma.$$
	Therefore, using Theorem \ref{maximumprinciple2} , we have for $a.e.$ $\mathbf{x}\in\Omega$ and for $t\in[0,T)$,
	$
	u_1(t,\mathbf{x})\leq v(t,\mathbf{x})
	$
	and
	$
	u_2(t,\mathbf{x})\leq z(t,\mathbf{x}).
	$
	Hence, to establish global existence, it is enough to prove that $u_3\in L^p(Q_{T})$ for sufficiently large $p$.

	Now, let $q>1$. Then,  let us multiply the equation for $i=3$ of System~\eqref{SystemChemReac-DiffFR2}  by $(u_3(t,\mathbf{x}))^{q}$, and then integrate over $Q_T$. We obtain
	\begin{equation}\label{xuq+intQT}
		\left\{ \begin{array}{ll} \dis
			\frac{1}{q+1} \int_{\Omega} (u_3(T,\mathbf{x}))^{q+1}d\mathbf{x}+ d_3 \iint_{Q_T} \big[(-\Delta)_{Sp}^{s_3} u_3 (t,\mathbf{x}) \big] (u_3(t,\mathbf{x}))^{q} d\mathbf{x}dt\\\\\dis\hskip2cm+\gamma\iint_{Q_T} (u_3(t,\mathbf{x}))^{q+\gamma} d\mathbf{x}dt
			\\\\\dis=\gamma\iint_{Q_T} (u_1(t,\mathbf{x}))^{\alpha} (u_2(t,\mathbf{x}))^{\beta} (u_3(t,\mathbf{x}))^q  d\mathbf{x}dt+\frac{1}{q+1} \int_{\Omega} (u_{03}(\mathbf{x}))^{q+1} d\mathbf{x}.
		\end{array}\right.
	\end{equation}
	Thanks to Stroock-Varopoulos inequality (see \cite[Lemma 6.2]{BonIbaIsp2023}), we have
	\begin{equation}
			\left\{
		\begin{array}{ll}\dis
		\iint_{Q_T} \big[(-\Delta)_{Sp}^{s_3} u_3 (t,\mathbf{x}) \big] (u_3(t,\mathbf{x}))^{q} d\mathbf{x}dt\\\\\dis \geq \frac{4q}{(q+1)^2}   \iint_{Q_T} \left| (-\Delta)_{Sp}^{\frac{s_3}{2}} (u_3(t,\mathbf{x}))^{\frac{q+1}{2}}  \right|^2 d\mathbf{x}dt\geq 0.
			\end{array}
		\right.
	\end{equation}
	In addition, using H\"{o}lder's inequality, we get
	\begin{equation}
			\left\{
		\begin{array}{ll}\dis
		\iint_{Q_T} (u_1(t,\mathbf{x}))^{\alpha} (u_2(t,\mathbf{x}))^{\beta} (u_3(t,\mathbf{x}))^q  d\mathbf{x}dt \\\\\dis\leq
		\|u_1\|_{L^{\alpha \sigma_1}\left(Q_T\right)}^{\alpha}\|u_2\|_{L^{\beta \sigma_2}\left(Q_T\right)}^{\beta}\|u_3\|_{L^{\gamma+q}\left(Q_T\right)}^q,	\end{array}
		\right.		
	\end{equation}
	where
	$$
	\frac{1}{\sigma_1}+\frac{1}{\sigma_2}+\frac{q}{q+\gamma}=1.
	$$
	Then,  we may choose $(\sigma_1,\sigma_2)$ such that $\sigma_1 \alpha \leq q+\gamma$  and  $\sigma_2\beta \leq q+\gamma$. This choice is valid, as by hypothesis $\gamma>\alpha+\beta$. Therefore, there exists $C_1>0$ such that
	\begin{equation}\label{Linclusion}
			\left\{
		\begin{array}{ll}\dis
		\iint_{Q_T} (u_1(t,\mathbf{x}))^{\alpha} (u_2(t,\mathbf{x}))^{\beta} (u_3(t,\mathbf{x}))^q  d\mathbf{x}dt \\\\\dis\leq C_1
		\|u_1\|_{L^{\gamma+q}\left(Q_T\right)}^{\alpha}\|u_2\|_{L^{\gamma+q}\left(Q_T\right)}^{\beta}\|u_3\|_{L^{\gamma+q}\left(Q_T\right)}^q.
			\end{array}
		\right.
	\end{equation}
	Since $s_3\leq\min\{s_1,s_2\}$, in accordance with Lemma \ref{FractionalDualityThm}, we get
	\begin{equation}\label{u1norm}
		\|u_1\|_{L^{\gamma+q}\left(Q_T\right)} \leq C_2\left(1+\|u_3\|_{L^{\gamma+q}\left(Q_T\right)}\right)
	\end{equation}
	and 
	\begin{equation} \label{u2norm}
		\|u_2\|_{L^{\gamma+q}\left(Q_T\right)} \leq C_3\left(1+\|u_3\|_{L^{\gamma+q}\left(Q_T\right)}\right) .
	\end{equation}
	Now, using the last two estimates \eqref{u1norm} and \eqref{u2norm}, we can rewrite \eqref{Linclusion} as
	\begin{equation}
		\dis
		\iint_{Q_T} (u_1(t,\mathbf{x}))^{\alpha} (u_2(t,\mathbf{x}))^{\beta} (u_3(t,\mathbf{x}))^q  d\mathbf{x}dt\leq C_4\left(1+\|u_3\|_{L^{\gamma+q}\left(Q_T\right)}\right)^{q+\alpha+\beta} .
	\end{equation}
	If $\|u_3\|_{L^{\gamma+q}\left(Q_T\right)} \leq 1$, the proof is complete. If not, further calculations are required. Indeed, we have 
	\begin{equation}
		\iint_{Q_T} (u_1(t,\mathbf{x}))^{\alpha} (u_2(t,\mathbf{x}))^{\beta} (u_3(t,\mathbf{x}))^q  d\mathbf{x}dt \leq C_5\|u_3\|_{L^{\gamma+q}\left(Q_T\right)}^{q+\alpha+\beta}.
	\end{equation}
	Then, thanks to \eqref{xuq+intQT}, we get 
	\begin{equation}
		\iint_{Q_T} (u_3(t,\mathbf{x}))^{q+\gamma}d\mathbf{x}dt\leq C_6\|u_3\|_{L^{\gamma+q}\left(Q_T\right)}^{q+\alpha+\beta}+\frac{1}{q+1} \int_{\Omega} (u_{03}(\mathbf{x}))^{q+1} d\mathbf{x}.
	\end{equation}
	In the sequel, we will apply the Young's inequality
	$ab\leq \ \varepsilon a^\eta+C(\varepsilon) b^\mu$, where $\frac{1}{\eta}+\frac{1}{\mu}=1$ and $\varepsilon\in(0,1)$. As $q+\gamma>q+\alpha+\beta$, then by choosing $\lambda:= \frac{q+\gamma}{q+\alpha+\beta}$, Young's inequality yields
	\begin{equation}
		(1-\varepsilon)	\times \iint_{Q_T} (u_3(t,\mathbf{x}))^{q+\gamma}d\mathbf{x}dt\leq \frac{1}{q+1} \int_{\Omega} (u_{03}(\mathbf{x}))^{q+1} d\mathbf{x}+C_7.
	\end{equation}
	Thus, we get
	\begin{equation}
		\|u_3\|_{L^{\gamma+q}\left(Q_T\right)}\leq C_8.
	\end{equation}
	Let us choose $\dfrac{\gamma+q}{\gamma}>\max\left\{\dfrac{N+2s_1}{2s_1},\dfrac{N+2s_2}{2s_2}\right\}$. Returning to Problems \eqref{P1} and \eqref{P2}, by applying Theorem \ref{BoundedSol} we obtain 
	\begin{equation}\label{u1bounded}
		\|u_1\|_{L^\infty(Q_T)}\leq C_9
	\end{equation}
	and 	\begin{equation}\label{u2bounded}
		\|u_2\|_{L^\infty(Q_T)}\leq C_{10}.
	\end{equation}
	Now, let us go back to the equation for $i=3$ of System \eqref{SystemChemReac-DiffFR2}. Thanks to the estimates \eqref{u1bounded} and \eqref{u2bounded}, we conclude that
	\begin{equation}
		\|u_3\|_{L^{\infty}\left(Q_T\right)}\leq C_{11}.
	\end{equation}
	Hence, we get $T_{\max}=+\infty$.
	\bigbreak      
	\noindent {\bf(ii)} Let  $s_3\geq\max\{s_1,s_2\} $ and $(\alpha, \beta,\gamma)\in  [1,+\infty)^2\times \{1\}$. Going back to Problem \eqref{P1}, we have for any $(t,\mathbf{x})\in Q_T$,
	\begin{equation}
		v(t,\mathbf{x})= T_{s_1}(t)u_{01}(\mathbf{x})+\alpha\int_0^tT_{s_1}(t-\tau)u_3(\tau,\mathbf{x})d\tau.
	\end{equation}
	By taking the $L^p$-norm and using the fact that $\{T_{s_1}(t)\}_{t\geq 0}$ is a contraction semigroup on $L^p(\Omega)$, we obtain
	\begin{equation}
		\left\{
		\begin{array}{ll}
			\dis	\|v(t,\cdot)\|_{L^p(\Omega)}&\leq\dis \left\|T_{s_1}(t)u_{01}\right\|_{L^p(\Omega)}+\alpha\int_0^t\|T_{s_1}(t-\tau)u_3(\tau,\cdot)\|_{L^p(\Omega)} d\tau\\\\\dis &\dis\leq\left\|u_{01}\right\|_{L^p(\Omega)}+\alpha\int_0^t\|u_3(\tau,\cdot)\|_{L^p(\Omega)} d\tau.
		\end{array}
		\right.
	\end{equation}
	As for $a.e.$ $\mathbf{x}\in\Omega$ and for any $t\geq 0$, we have
	$
	u_1(t,\mathbf{x})\leq v(t,\mathbf{x})
	$, we deduce that
	\begin{equation}\label{inequ_1}
		\|u_1(t,\cdot)\|_{L^p(\Omega)} \leq\left\|u_{01}\right\|_{L^p(\Omega)}+\alpha\int_0^t\|u_3(\tau,\cdot)\|_{L^p(\Omega)} d\tau.
	\end{equation}
	By applying H\"{o}lder's inequality, we have
	\begin{equation}
		\int_0^t\|u_3(\tau,\cdot)\|_{L^p(\Omega)} d\tau \leq t^{\frac{1}{p^\prime}} \left(\int_0^t\int_\Omega |u_3(\tau,\mathbf{x})|^p d\tau d\mathbf{x} \right)^{\frac{1}{p}},
	\end{equation}
	where $p^\prime=\dfrac{p}{p-1}$.
	Moreover, by using Lemma \ref{FractionalDualityThm} and the fact that $s_3\geq s_1$, we get
	\begin{equation}
		\int_0^t\|u_3(\tau,\cdot)\|_{L^p(\Omega)} d\tau \leq C_{12}t^{\frac{1}{p^\prime}}  \left[1+\left(\int_0^t \int_{\Omega} |u_1(\tau,\mathbf{x})|^p d\tau d \mathbf{x}\right)^{\frac{1}{p}}\right].
	\end{equation}
	Now, for any $t \in(0, T]$, let us denote $g(t):=\dis\int_{\Omega}|u_1(t, \mathbf{x})|^p d\mathbf{x}$. Thus, we can rewrite \eqref{inequ_1} as
	\begin{equation}\label{ineqh}
		g(t)^{\frac{1}{p}} \leq C_{13}\left[1+\left(\int_0^t g(\tau) d\tau\right)^{\frac{1}{p}}\right].
	\end{equation}
	This yields
	\begin{equation}
		g(t) \leq C_{14}\left[1+\int_0^t g(\tau) d\tau \right].
	\end{equation}
	By Gronwall's lemma, we deduce that
	\begin{equation}\label{estu1}
		\|u_1\|_{L^p\left(Q_T\right)}\leq C_{15} .
	\end{equation}
	By applying the same reasoning to $u_2$ in place of $u_1$ and using the fact that $s_3\geq s_2$, we derive that
	\begin{equation} \label{estu2}
		\|u_2\|_{L^p\left(Q_T\right)}\leq C_{16}.
	\end{equation}
	Thanks to the estimates \eqref{estu1} and \eqref{estu2}, we may choose $q>\frac{N+2s_3}{2s_3}$ such that
	\begin{equation}
		\|u_1^{\alpha} u_2^{\beta}\|_{L^q\left(Q_T\right)} \leq C_{17}.
	\end{equation}
	Now, let us go back to the equation for $i=3$ of System \eqref{SystemChemReac-DiffFR2}. By using Theorem \ref{BoundedSol}, we obtain
	\begin{equation}
		\|u_3\|_{L^{\infty}\left(Q_T\right)}\leq C_{18}.
	\end{equation}
	Finally $T_{\max}=+\infty$, which concludes the proof.

\end{proof}

	\subsection{Global existence : case of triangular structure}\label{triangular-structure}

Our second result on global existence pertains to System
\eqref{ReactiondiffusionSystem3} with $m\geq 2$, under an additional assumption on the source terms; referred to as the {\it triangular structure}.

\begin{theorem}\label{MainTh2} % Theorem 4.2 {Extension2}
	For each $i\in\llbracket1,m\rrbracket$, with $m\geq2$, let $u_{0i} \in L^\infty(\Omega)^+$. Moreover, assume that $s_i\in (0,1)$, where $\forall\ell\in\llbracket1,m-1\rrbracket$, $s_\ell\leq s_{\ell+1}$. In addition to {\bf (P)},  let us assume that the $f_i$'s are  
	at most polynomial satisfying the following so-called  ``triangular structure":
	\begin{equation}\label{TriangStructure}\tag{TS}
		\left\{
		\begin{array}{lll}
			\text{There exist a vector} \; \mathbf{b}\in\mathbb{R}_+^m\; \text{and a lower triangular invertible matrix}\; \\Q\in {\mathcal{M}_m(\mathbb{R}_+)}\; (q_{ii}\neq 0), \text{such that} \;\;	\forall\mathbf{r}=(r_1,\ldots,r_m)\in \mathbb{R}_+^m,\\Q\mathbf{f}(\mathbf{r})\leq \Big[1+\dis\sum\limits_{i=1}^m r_i\Big]\mathbf{b}\;\;\text{where}\;\; \mathbf{f}(\mathbf{r})=(f_1(\mathbf{r}),\ldots,f_m(\mathbf{r}))^T.
		\end{array}
		\right.
	\end{equation}

	Then, System \eqref{ReactiondiffusionSystem3} admits a unique nonnegative global strong  solution.
\end{theorem}

\begin{proof}
	Let us denote
	$Q=(q_{ij})_{i,j=1,\ldots,m}$ and	$\mathbf{b}=(b_1,\ldots,b_m)^T$.
	Thanks to Lemma \ref{LocalExistence}, System \eqref{ReactiondiffusionSystem3} has a unique strong solution $(u_1,\cdots,u_m)$ on $Q_{T_{\max}}$.  In the subsequent, we will prove that this solution is global, meaning that $T_{\max}=+\infty$.
	Let $T\in(0,T_{\max})$ and $t\in (0,T]$.\\
	\noindent	For each $i\in\llbracket1,m\rrbracket$, let $v_i$ be the weak solution to
	\begin{equation}\tag{$P_i$}
		\label{Testfunc}
		\left\{  
		\begin{array}{rcll}
			\partial_t v_i(t,\mathbf{x})+d_i (-\Delta)_{Sp}^{s_i} v_i(t,\mathbf{x}) &=& \Big[1+\dis\sum\limits_{j=1}^m u_j(t,\mathbf{x})\Big] \dfrac{b_i}{q_{ii}}, &(t,\mathbf{x})\in Q_T,\\
			\mathcal{B}[v_i(t,\mathbf{x})]&=&0, &(t,\mathbf{x})\in \Sigma_T,\\
			v_i(0,\mathbf{x})&=&0,& \mathbf{x}\in\Omega.
		\end{array}
		\right.
	\end{equation}	
	Therefore, we have
	\begin{equation}
		\label{ineqlasttheo5}
		v_i(t,\mathbf{x})=\dis\int_0^t T_{s_i}(t-\tau)\Big[1+  \dis\sum\limits_{j=1}^mu_j(\tau,\mathbf{x}) \Big]\dfrac{b_i}{q_{ii}}d\tau,\quad (t,\mathbf{x})\in [0,T)\times\Omega.
	\end{equation}
	\\
	Given \eqref{TriangStructure}, the following inequality holds:
	\begin{equation}\label{sum-}
		\sum\limits_{j=1}^i q_{ij} f_j(u_1(t,\textbf{x}),\ldots,u_m(t,\textbf{x}))\leq \Big[ 1+    \sum\limits_{j=1}^m  u_j (t,\textbf{x})  \Big] b_i.
	\end{equation}
	---	Let us multiply both the first equation of $(P_1)$ and the equation for $i=1$ of System \eqref{ReactiondiffusionSystem1} by $q_{11}$, then subtract the former from the latter. We obtain
	\begin{equation}\label{ineq1}
		q_{11}\Big[\partial_\tau(u_1-v_1)(\tau,\textbf{x})+d_1(-\Delta)_{Sp}^{s_1}(u_1-v_1)(\tau,\textbf{x})  \Big]\leq 0,
	\end{equation}	
	after replacing $t$ by $\tau$.
	In light of the duality concept outlined in Lemma \ref{FractionalDualityThm},
	let $\phi$ be a nonnegative regular function and let $\mathcal{V}(\tau,\textbf{x})$ be the solution to {Problem $(P_{\phi,t})$} with $s=s_1$ and $d={d_1}$. Subsequently, we multiply inequality  \eqref{ineq1} by $\mathcal{V}$ and integrate over $Q_t$. As $q_{11}>0$, we get
	\begin{equation}\label{ineq2}
		\iint_{Q_t}\partial_t(u_1-v_1)(\tau,\textbf{x})\mathcal{V}(\tau,\textbf{x})d\tau d\textbf{x}+\iint_{Q_t}d_1(-\Delta)_{Sp}^{s_1}(u_1-v_1)(\tau,\textbf{x})\mathcal{V}(\tau,\textbf{x})d\tau d\textbf{x} \leq 0.
	\end{equation}	
	Then,
	\begin{equation}\label{ineq3}
			\left\{
		\begin{array}{ll}\dis
		\iint_{Q_t}-(u_1-v_1)(\tau,\textbf{x})\partial_\tau\mathcal{V}(\tau,\textbf{x})d\tau d\textbf{x}\\\\\dis\hskip2cm+\iint_{Q_t}d_1(u_1-v_1)(\tau,\textbf{x})(-\Delta)_{Sp}^{s_1}\mathcal{V}(\tau,\textbf{x})d\tau d\textbf{x}\\\\\dis \leq \int_{\Omega} u_{01}(\textbf{x}) \mathcal{V}_0 (\textbf{x})d\textbf{x}.
			\end{array}
		\right.
	\end{equation}	
	Therefore, we have
	\begin{equation}\label{ineq4}
		\iint_{Q_t}(u_1-v_1)(\tau,\textbf{x})\phi(\tau,\textbf{x})d\tau d\textbf{x}\leq \int_{\Omega} u_{01}(\textbf{x}) \mathcal{V}_0 (\textbf{x})d\textbf{x}.
	\end{equation}	
	{Consider a sufficiently large $p>1$ and let $p^\prime=\dfrac{p}{p-1}$.} By applying H\"{o}lder's inequality and the estimate \eqref{DualityInequality},  we obtain
	\begin{equation}\label{ineq5}
		%		\left\{
		%		\begin{array}{lll}
			%	\dis	
			\iint_{Q_t}(u_1-v_1)(\tau,\textbf{x})\phi(\tau,\textbf{x})d\tau d\textbf{x}\leq\dis \|u_{01}\|_{L^p(\Omega)} \|\mathcal{V}_0\|_{L^{p^\prime}(\Omega)}
			\leq  C_{19}  \|\phi\|_{L^{p^\prime}(Q_t)}.
			%		\end{array}
		%		\right.
	\end{equation}	
	By duality, it follows  that
	\begin{equation}
		\label{ineq6}
		\|(u_1-v_1)^+\|_{L^p(Q_t)}\leq C_{19} .
	\end{equation}
	Consequently, we obtain
	\begin{equation}\label{ineq7}
		\|u_1\|_{L^p(Q_{t})} \leq 
		\|v_1\|_{L^p(Q_{t})}
		+	\|(u_1-v_1)^+\|_{L^p(Q_{t})}\leq
		C_{20} [ 1+ \|	v_1 \|_{L^p(Q_{t})} ].
	\end{equation}
	\\
	--- Now, consider $i\geq 2$.	
	By multiplying the first equation of  \eqref{Testfunc} and the equation for $i$ of System \eqref{ReactiondiffusionSystem3} by $q_{ii}$, and then subtracting the former from the latter, we get

	\begin{equation}
		\label{ineqlasttheo1}
			\left\{
		\begin{array}{ll}\dis
		q_{ii}\Big[\partial_\tau(u_i-v_i)(\tau,\textbf{x})+d_i(-\Delta)_{Sp}^{s_i}(u_i-v_i)(\tau,\textbf{x})  \Big]\\\\\dis\leq -\dis\sum\limits_{j=1}^{i-1} q_{ij}\Big[ \partial_\tau u_j (\tau,\textbf{x}) +d_j(-\Delta)_{Sp}^{s_j} u_j(\tau,\textbf{x})  \Big],
			\end{array}
		\right.
	\end{equation}
	after replacing $t$ by $\tau$.
	As above,
	let $\phi$ be a nonnegative regular function and $\mathcal{V}(\tau,\textbf{x})$ be the solution to {Problem $(P_{\phi,t})$} with $s=s_i$ and $d={d_i}$. Next, we multiply both sides of inequality \eqref{ineqlasttheo1} by $\mathcal{V}$ and then integrate over $Q_t$, resulting in %for any $t\in (0,T_{\max})$. 
	\begin{equation}\label{I1I2} \left\{\begin{array}{ll}
			%\underbrace{
				q_{ii}	\dis\iint_{Q_t}\partial_\tau (u_i-v_i) (\tau,\textbf{x}) \mathcal{V}(\tau,\textbf{x})d\tau d\textbf{x}\\\\\dis\hskip2cm+q_{ii}\iint_{Q_t}d_i(-\Delta)_{Sp}^{s_i}(u_i-v_i)(\tau,\textbf{x}) \mathcal{V}(\tau,\textbf{x})d\tau d\textbf{x}
				%}_{I_1}
			\\\\\leq-\dis \sum\limits_{j=1 }^{i-1} q_{ij}
			%\underbrace{
				\iint_{Q_t}\dis	\partial_\tau u_j (\tau,\textbf{x})\mathcal{V}(\tau,\textbf{x})d\tau d\textbf{x}\\\\\dis\hskip2cm-\sum\limits_{j=1 }^{i-1} q_{ij}d_j\iint_{Q_t}(-\Delta)_{Sp}^{s_j}u_j(\tau,\textbf{x})\mathcal{V}(\tau,\textbf{x})d\tau d\textbf{x}.
				%}_{I_2}.
		\end{array}\right.
	\end{equation}

	Thus, we get
	\begin{equation}
		\left\{
		\begin{array}{ll}	
			\displaystyle	q_{ii} \iint_{Q_t}(u_i-v_i)(\tau,\textbf{x})\phi(\tau,\textbf{x}) d\tau d\textbf{x}\\\leq \dis \sum\limits_{j=1 }^{i-1} q_{ij}\int_{\Omega} u_{0j} (\textbf{x})\mathcal{V}_{0}(\textbf{x}) d\textbf{x}\\\\\dis\hskip2cm\dis+\sum\limits_{j=1}^{i-1} q_{ij} \iint_{Q_t}u_j(\tau,\textbf{x})\left( \partial_t\mathcal{V}(\tau,\textbf{x})-d_j(-\Delta)_{Sp}^{s_j}\mathcal{V}(\tau,\textbf{x}) \right)d\tau d\textbf{x}.
		\end{array} 
		\right.
	\end{equation}
	Consistent with the above, let us choose $p>1$ sufficiently large. %and set 
	% $p^\prime=\dfrac{p}{p-1}$.
	Applying H\"{o}lder's inequality along with the estimate \eqref{DualityInequality} yields
	\begin{equation}
		\left\{
		\begin{array}{ll}
			\displaystyle	\iint_{Q_t}(u_i-v_i)(\tau,\textbf{x})\phi(\tau,\textbf{x}) d\tau d\textbf{x}\\\\\displaystyle\leq C_{21}  \dis \sum\limits_{j=1}^{i-1} \|u_{0j}\|_{L^p(\Omega)}\|\mathcal{V}_0\|_{L^{p^\prime}(\Omega)}\\\\\dis\hskip2cm+C_{21}   \dis \sum\limits_{j=1}^{i-1}\|u_j\|_{L^p(Q_t)}\left(\left\| \partial_\tau\mathcal{V} \right\|_{L^{p^\prime}(Q_t)} + d_j \left\| (-\Delta)_{Sp}^{s_j}\mathcal{V} \right\|_{L^{p^\prime}(Q_t)}  \right)
			\\\\\displaystyle\leq C_{22}   \dis \sum\limits_{j=1}^{i-1} \|u_{0j}\|_{L^p(\Omega)}\|\phi\|_{L^{p^\prime}(Q_t)}\\\\\dis\hskip2cm+C_{22}   \dis \sum\limits_{j=1}^{i-1}\|u_j\|_{L^p(Q_t)}\left(\left\| \phi \right\|_{L^{p^\prime}(Q_t)} + d_j \left\| (-\Delta)_{Sp}^{s_j}\mathcal{V} \right\|_{L^{p^\prime}(Q_t)}  \right).
		\end{array}
		\right.
	\end{equation}	
	By assumption, we have $0< s_1\leq \cdots \leq s_m<1$. Then, for any fixed $i\in\llbracket2,m\rrbracket$ and any $j\in\llbracket1, i-1\rrbracket$, it follows that $s_j\leq s_i$. Therefore, thanks to Proposition \ref{estimation}, for any $j\in\llbracket1, i-1\rrbracket$,
	\begin{equation}
			\left\{
		\begin{array}{ll}
			\displaystyle
		\left\| (-\Delta)_{Sp}^{s_j}\mathcal{V} \right\|_{L^{p^\prime}(Q_t)} & \leq  C \|(-\Delta)_{Sp}^{s_i}\mathcal{V}\|_{L^{p^\prime}(Q_t)}^{\frac{s_j}{s_i}}  \|\mathcal{V}\|_{L^{p^\prime}(Q_t)}^{\frac{s_i-s_j}{s_i}}\\\\&
		\displaystyle\leq C \|\phi\|^{\frac{s_j}{s_i}} _{L^{p^\prime}(Q_t)} \|\phi\|^{\frac{s_i-s_j}{s_i}} _{L^{p^\prime}(Q_t)}=C \left\| \phi \right\|_{L^{p^\prime}(Q_t)}, 
			\end{array}
		\right.
	\end{equation}
	where the second inequality follows from the estimate \eqref{DualityInequality}. Consequently, we obtain
	\begin{equation}
		\left\{
		\begin{array}{ll}
			\displaystyle	\iint_{Q_t}(u_i-v_i)(\tau,\textbf{x})\phi(\tau,\textbf{x}) d\tau d\textbf{x}
			\\\\ \dis	\displaystyle\leq \displaystyle C_{22} \dis \sum\limits_{j=1 }^{i-1}\|u_{0j}\|_{L^p(\Omega)} \|\phi\|_{L^{p^\prime}(Q_t)} +C_{23}  \dis \sum\limits_{j=1 }^{i-1} \|u_j\|_{L^p(Q_t)}\|\phi\|_{L^{p^\prime}(Q_t)}\\\\
			\displaystyle\leq C_{24} \Big(1+\dis \sum\limits_{j=1}^{i-1} \|u_j\|_{L^p(Q_t)}\Big)\|\phi\|_{L^{p^\prime}(Q_t)}.
		\end{array}
		\right.
	\end{equation}		
	Since the function $\phi$ is nonnegative and regular, we deduce by duality that %for any $p>1$
	\begin{equation}
		\label{ineqlasttheo2}
		\|(u_i-v_i)^+\|_{L^p(Q_t)}\leq C_{24}  \Big( 1+ \dis\sum\limits_{j=1}^{i-1}\|u_j\|_{L^p(Q_t)}\Big).
	\end{equation}
	Therefore, given $u_i\leq v_i + (u_i-v_i)^+$ and applying elementary induction on $i$, we get %for each  $i=1,\ldots,m$,
	\begin{equation}
		\label{ineqlasttheo2-}
		\|u_i\|_{L^p(Q_t)}\leq C_{25}  \Big( 1+ \dis\sum\limits_{j=1 }^i\|v_j\|_{L^p(Q_t)}\Big),
	\end{equation}
	for any $i\in\llbracket2,m\rrbracket$. Thanks to \eqref{ineq7}, inequality  \eqref{ineqlasttheo2-} is valid 	for any $i\in\llbracket1,m\rrbracket$.
	Hence, summing over $i$ and raising to the $p$th power yields
	\begin{equation}
		\label{ineqlasttheo3}
		\Big\|\dis\sum\limits_{i=1}^mu_i\Big\|^p_{L^p(Q_t)}\leq C_{26}  \Big( 1+ \sum\limits_{i=1}^m\|v_i\|^p_{L^p(Q_t)}\Big).
	\end{equation}
	\noindent 
	By taking the \(L^p(\Omega)\)-norm of \eqref{ineqlasttheo5} and using the estimate \eqref{Ultracontractive} with \(q = p\) together with H\"{o}lder's inequality, we obtain
	
	\begin{equation}
		\label{ineqlasttheo6}
		\|v_i(t,\cdot)\|^p_{L^p(\Omega)}\leq C_{27} \Big[ 1+   \Big\|\dis\sum\limits_{i=1}^mu_i\Big\|^p_{L^p(Q_t)}   \Big].
	\end{equation}
	Then, by using \eqref{ineqlasttheo3}, we derive that
	\begin{equation}
		%		\label{ineqlasttheo7}
		%		\left\{
		%		\begin{array}{lll}
			\dis\sum\limits_{i=1}^m\|v_i(t,\cdot)\|^p_{L^p(\Omega)}\leq C_{28}  \left[ 1+  	\dis  \int_0^t\sum\limits_{i=1}^m \|v_i(\tau,\cdot)\|^p_{L^p(\Omega)}   d\tau \right].
			%			\end{array}
		%		\right.
	\end{equation}	
	Hence, for any  $i\in\llbracket1,m\rrbracket$,  $v_i(t, \cdot)$ is bounded in $L^p(\Omega)$ for any $p<+\infty$. Returning to \eqref{ineqlasttheo3}, $\dis\sum\limits_{i=1}^m u_i(t, \cdot)$ is also bounded in $L^p(\Omega)$ for any $p<+\infty$, and consequently, so is $u_i(t, \cdot)$.
	\smallbreak\noindent
	Since the $f_i$'s are at most polynomial, it follows that $f_i \in L^q\left(Q_T\right)$ for sufficiently large $q$.
	Therefore, we can select
	$q>\max\limits_{i\in\llbracket1,m\rrbracket}\Big\{\dfrac{N+2s_i}{2s_i}\Big\}$ to ensure, by Theorem \ref{BoundedSol}, that for any $i\in\llbracket1,m\rrbracket$
	\begin{equation}
		\|u_i\|_{L^\infty(Q_{T})}\leq C_{29} .
		%\quad\forall i=1,\ldots,m.
	\end{equation}
	Finally, $T_{\max}=+\infty$.
\end{proof}

\section{Numerical simulations for System $(S_{\alpha,\beta,\gamma})$}

The purpose of this section is to numerically investigate the three-species
fractional reaction–diffusion system $(S_{\alpha,\beta,\gamma})$, introduced in
Subsection~\ref{reversible-chemical-reaction}, in a parameter regime that remains
open from a theoretical point of view.  
For convenience, we recall that the system is given by
\begin{equation}\label{System-FR-Num}\tag{$S_{\alpha,\beta,\gamma}$}
	\left\{
	\begin{array}{rcll}
		\partial_t u_1 + d_1(-\Delta)_{Sp}^{s_1} u_1 &=& \alpha\,
		\bigl(u_3^{\gamma}-u_1^{\alpha}u_2^{\beta}\bigr),
		& (t,\mathbf{x})\in (0,T)\times \Omega,\\[3pt]
		\partial_t u_2 + d_2(-\Delta)_{Sp}^{s_2} u_2 &=& \beta\,
		\bigl(u_3^{\gamma}-u_1^{\alpha}u_2^{\beta}\bigr),
		& (t,\mathbf{x})\in (0,T)\times \Omega,\\[3pt]
		\partial_t u_3 + d_3(-\Delta)_{Sp}^{s_3} u_3 &=& -\gamma\,
		\bigl(u_3^{\gamma}-u_1^{\alpha}u_2^{\beta}\bigr),
		& (t,\mathbf{x})\in (0,T)\times \Omega,\\[3pt]
		\partial_{\nu} u_i &=& 0, & (t,\mathbf{x})\in (0,T)\times \partial\Omega,\ i=1,2,3,\\[3pt]
		u_i(0,\mathbf{x}) &=& u_{0i}(\mathbf{x}), & \mathbf{x}\in\Omega,
	\end{array}
	\right.
\end{equation}
where $0<s_i<1$, $d_i>0$, and $(\alpha,\beta,\gamma)\in[1,+\infty)^3$. As mentioned earlier, the reaction terms satisfy the structural assumptions
$(\mathbf P)$ and $(\mathbf M)$.

\medskip

Let us recall that,
in the classical diffusion setting ($s_i=1$), it is known that 
any global solution (when it exists) converges in $L^1$ towards an
equilibrium state as time goes to infinity. More precisely, let us denote
\[
\overline{u_{0i}}
= \frac{1}{|\Omega|}
\int_\Omega u_{0i}(\mathbf{x})\,d\mathbf{x},
\qquad i=1,2,3,
\]
and define the two conserved
mean quantities
\[
M_1 := \gamma\,\overline{u_{01}} + \alpha\,\overline{u_{03}}, 
\qquad
M_2 := \gamma\,\overline{u_{02}} + \beta\,\overline{u_{03}}.
\]
Then the equilibrium $(u_{1,\infty},u_{2,\infty},u_{3,\infty})$ is
uniquely determined as the solution of the algebraic system
\begin{equation}\label{Eq-Equilibrium-Pierre}
	\left\{
	\begin{aligned}
		\gamma u_{1,\infty} + \alpha u_{3,\infty} &= M_1,\\
		\gamma u_{2,\infty} + \beta  u_{3,\infty} &= M_2,\\
		u_{1,\infty}^{\alpha}\,u_{2,\infty}^{\beta} &= u_{3,\infty}^{\gamma}.
	\end{aligned}
	\right.
\end{equation}
From the first two linear equations we explicitly obtain 
\(u_{1,\infty}\) and \(u_{2,\infty}\) as functions of 
\(u_{3,\infty}\). Inserting these formulas into the third 
equation yields a nonlinear equation involving only 
\(u_{3,\infty}\). This reduces the problem to a single scalar
equation for $w:=u_{3,\infty}$,
\[
\Bigl(\tfrac{M_1-\alpha w}{\gamma}\Bigr)^{\alpha}
\Bigl(\tfrac{M_2-\beta w}{\gamma}\Bigr)^{\beta}
= w^{\gamma},
\qquad
0 < w < 
\min\Bigl\{\tfrac{M_1}{\alpha},\,\tfrac{M_2}{\beta}\Bigr\},
\]
which admits a unique solution.
The equilibrium values are then recovered explicitly as
\[
u_{1,\infty} = \frac{M_1-\alpha w}{\gamma},
\qquad
u_{2,\infty} = \frac{M_2-\beta w}{\gamma}.
\]
In particular, the strong solution $(u_1,u_2,u_3)$ satisfies
\[
u_i(t,\cdot)\xrightarrow[t\to+\infty]{} u_{i,\infty}
\quad \text{in } L^1(\Omega).
\]
For complete proofs of these statements, we refer to  \cite{PierSuzUma2018,FellLaam2016}.

\medskip

To the best of our knowledge, no analogous convergence
result is known when the diffusion is governed by spectral fractional Laplacians $(-\Delta)_{Sp}^{s_i}$. In Theorem~\ref{MainTh1}, we established the global existence of strong 
solutions to System \eqref{System-FR-Num} whenever
$
s_3 \le \min\{s_1,s_2\}$ and
$
\gamma > \alpha+\beta$.
This may be viewed as a fractional extension of the classical 
result obtained in \cite[Theorem 1]{LaamActa2011} in the case of the classical Laplacian 
($s_i=1$), where the same condition $\gamma>\alpha+\beta$ ensures 
global existence of strong solutions.
In contrast, the parameter regimes
$
s_3 > \min\{s_1,s_2\}$
and 
$ \gamma\leq\alpha+\beta$,
remain open from an analytical point of view.
These are precisely the configurations we aim to explore numerically.

\medbreak
Nonetheless, it should be emphasized that the numerical approximation of such models is challenging and imposes significant computational constraints. To tackle these challenges, we adopt a Fourier spectral method introduced in \cite{BuenoKayBurr2014}, which provides an appealing and straightforward approach for solving fractional problems of type  \eqref{System-FR-Num} in bounded rectangular domains  of $\mathbb{R}^N$, for $N\in\llbracket1,3\rrbracket$. Indeed, their proposed schemes yield a fully diagonal representation of the $SFL$, enhancing both accuracy and computational efficiency.

\subsection{Description of the numerical method in 3D}

Now, let us give a brief overview of the chosen numerical method in the case of the homogeneous Neumann boundary condition. Let us  underline that the method also applies when considering the homogeneous Dirichlet boundary condition instead of the Neumann one.

\medbreak
\noindent
{\bf $\bullet$ Spatial discretization.} 
Let $\Omega=(c_1,c_2)^3$ be a bounded cubic domain, and consider the following fractional diffusion equation:
\begin{equation}\label{HeatEq3D}
	\left\{
	\begin{array}{rclll}
		\partial_t w(t,\mathbf{x}) + d(-\Delta)_{Sp}^s w(t,\mathbf{x}) &=& 0, & (t,\mathbf{x})\in (0,T)\times\Omega,\\[3pt]
		\partial_{\nu} w(t,\mathbf{x}) &=& 0, & \mathbf{x}\in \partial\Omega,~t\in(0,T),\\[3pt]
		w(0,\mathbf{x}) &=& w_0(\mathbf{x}), & \mathbf{x}\in \Omega.
	\end{array}
	\right.
\end{equation} 
Therefore, the analytical solution to Problem \eqref{HeatEq3D} can be written as
\begin{equation}\label{eq-}
	w(t,\mathbf{x})
	=\sum_{\mathbf{k}\in\mathbb{N}^3} w_{\mathbf{k}}(t) e_{\mathbf{k}}(\mathbf{x})
	=\sum_{\mathbf{k}\in\mathbb{N}^3} 
	w_{\mathbf{k}}(0) e^{-d\,t\,\lambda_{\mathbf{k}}^{s}} e_{\mathbf{k}}(\mathbf{x}), 
\end{equation}
where  $	\mathbf{k}=(k_1,k_2,k_3)\in \mathbb{N}^3,$
\[
\lambda_{\mathbf{k}} = 
\left(\frac{\pi}{c_2-c_1}\right)^2 (k_1^2+k_2^2+k_3^2)
\;\text{and} \;
e_{\mathbf{k}}(\mathbf{x}) = 
\prod_{j=1}^3 \sqrt{\frac{2}{c_2-c_1}}\,
\cos\!\left(\frac{k_j\pi(x_j-c_1)}{c_2-c_1}\right).
\]
The Fourier spectral method under discussion involves approximating the series expansion of \eqref{eq-} by using a finite set of orthonormal trigonometric eigenfunctions $\left\{e_{\mathbf{k}}\right\}_{\mathbf{k}\in\mathbb{N}^3}$, where the number of these functions is equal to the number of discretization points. Then, we have
\begin{equation}\label{eq--}
	w(t,\mathbf{x}) \approx 
	\sum_{k_1=0}^{K}\sum_{k_2=0}^{K}\sum_{k_3=0}^{K}
	w_{\mathbf{k}}(0)\, e^{-d\,t\,\lambda_{\mathbf{k}}^{s}}\, e_{\mathbf{k}}(\mathbf{x}).
\end{equation}

\medbreak
The coefficients $w_{\mathbf{k}}(0)$ in \eqref{eq--}, along with the inverse reconstruction of $w$, 
can be efficiently computed using established algorithms such as the direct and inverse 
Discrete Sine/Cosine Transforms, depending on the boundary conditions; 
see, for instance, \cite{BrigHen1995,BuenoPer2006}. 
In \cite{BuenoKayBurr2014}, the authors have illustrated the simplicity of applying this method.
Additionally,  $\Omega$ is discretized using a uniform Cartesian grid with
mesh size $h_x := (c_2-c_1)/(K+1)$.
The grid points are given by
\[
\mathbf{x}=(x_i,y_j,z_k) = (c_1 + (i-1) h_x+h_x/2, c_1 + (j-1) h_x+h_x/2, c_1 + (k-1) h_x+h_x/2),
\]
for each $i,j,k \in \llbracket 1,K+1 \rrbracket$.

\smallbreak 
\textit{Comment.}
Problem \eqref{HeatEq3D} is recalled here because it contributes to the construction of the numerical scheme. Indeed, the operator $(-\Delta)_{Sp}^s$ is diagonal in the cosine basis, so each Fourier mode evolves independently and the action of the fractional Laplacian is computed exactly in spectral space. In the following steps, we add the reaction terms to this diffusive structure.

\medbreak
\noindent
{\bf $\bullet$ Time discretization.}  Let us go back to System \eqref{System-FR-Num}. Let us consider a discretization of time as $0=t_0<t_1<\cdots<t_{\widetilde{N}}=T$ with $h_t:=t_{n+1}-t_n$ for $n\in \llbracket0,\widetilde{N}-1\rrbracket$.
Following the approach in \cite{BuenoKayBurr2014}, we apply a backward Euler
discretization to the time derivative. Let us consider $u_i^n(\cdot)\approx u_i(t_n,\cdot)$ for $i=1,2,3$. In each time interval $\left[t_n, t_{n+1}\right]$, the nonlinear terms are handled using the following fixed point iteration: given $(u_1^n,u_2^n,u_3^n)$, initialize $(u_1^{n+1,0} ,u_2^{n+1,0},u_3^{n+1,0} )  :=  (u_1^n,u_2^n,u_3^n)$, and for each $\ell\in\llbracket1,L\rrbracket$ 
\begin{equation}\label{eq---}
	\left\{
	\begin{array}{lll}
\text{find }\; 
	(u_1^{n+1, \ell},u_2^{n+1, \ell},u_3^{n+1, \ell}) \; \text{ such that}
	\\\\
		\dfrac{u_i^{n+1, \ell}-u_i^n}{h_t}=-d_i(-\Delta)_{Sp}^{s_i} u_i^{n+1, \ell}+f_i\left(u_1^{n+1, \ell-1},u_2^{n+1, \ell-1},u_3^{n+1, \ell-1}\right)
	\end{array}
	\right.
\end{equation}
where $i=1,2,3$ and $L$ is to be determined later. Notably, $L=1$ corresponds to a fully explicit treatment of the nonlinear term, whereas a sufficiently large $L$ makes the method fully implicit.

Now, projecting both sides of \eqref{eq---} onto the eigenfunction basis  yields, for ${\mathbf{k}}\in\llbracket1,K\rrbracket^3$ and $i=1,2,3$,

$$
\frac{{u_i}_{\mathbf{k}}^{n+1, \ell}-{u_i}_{\mathbf{k}}^n}{h_t}=-d_i \lambda_{\mathbf{k}}^{s_i} {u_i}_k^{n+1, \ell}+{f_i}_{\mathbf{k}}\left(u_1^{n+1, \ell-1},u_2^{n+1, \ell-1},u_3^{n+1, \ell-1}\right)
$$
where ${f_{i}}_{\mathbf{k}}$ is the ${\mathbf{k}}$th Fourier coefficient of the $f_i$. Hence, after rearranging the terms, we obtain

$$
{u_i}_{\mathbf{k}}^{n+1, \ell}=\frac{1}{1+d_i \lambda_{\mathbf{k}}^{s_i}  h_t}\left[{u_i}_{\mathbf{k}}^n+h_t {f_{i}}_{\mathbf{k}}\left(u_1^{n+1, \ell-1},u_2^{n+1, \ell-1},u_3^{n+1, \ell-1}\right)\right].
$$

\subsection{Numerical results in 3D}
Now, using the chosen numerical scheme, let us showcase some numerical simulations. Our goal is to investigate whether a nonnegative solution to System \eqref{System-FR-Num} may
persist for sufficiently large times in the open theoretical regime where no global
existence result is currently known:
\begin{equation}\label{Hypo}
	s_3>\min\{s_1,s_2\},
	\qquad 
	\alpha+\beta\geq\gamma.
\end{equation}

\medbreak
To this end, we choose  $\Omega=(-5,5)^3$ and diffusion coefficients
$
d_1=\tfrac{3}{4}, \; d_2=\tfrac{1}{4}, \; d_3=\tfrac{1}{2}.
$
The initial data are taken as
\[
u_{0i}(x,y,z)
=
\Bigl[\frac12(1+\cos(\pi x))(1+\cos(\pi y))(1+\cos(\pi z))\Bigr]^{\,s_i},
\qquad i=1,2,3.
\]
Furthermore, we set
$
s_1=\tfrac{1}{2}, \; s_2=\tfrac{3}{4}, \; s_3=\tfrac{9}{10}$ and
$
\alpha=\beta=\gamma=2
$
so that the assumptions in \eqref{Hypo} are fully satisfied.

\medbreak In the numerical simulations below, we proceed as follows. First, we compute the solutions for intermediate times in order to observe their behaviour. % and to detect any possible instabilities.  
In order
to assess whether the solution tends towards a spatially homogeneous state, we simulate System \eqref{System-FR-Num} up to a large final time. 
Then, we compute the corresponding equilibrium 
\((u_{1,\infty},u_{2,\infty},u_{3,\infty})\)
by solving the system~\eqref{Eq-Equilibrium-Pierre}, exactly as in the case of the classical Laplacian. 
Although no analogous result is currently available for the fractional case, this procedure provides a natural reference state for comparison. Therefore, we compare the fractional numerical solutions to these constants. Finally, to test the robustness of our observations, we evaluate the stability of the numerical solution with respect to the discretization parameters (spatial resolution \(K\), number of iterations \(L\) and time step \(h_t\)), and monitor the \(L^2\)-norm of the solution at a large final time.

To visualize the dynamics of the system, we plot two-dimensional slices of the numerical solution on the mid-plane \(z=0\). Figures~\ref{fig:init}, \ref{fig:time}  and~\ref{fig:L1cv} correspond to some
numerical experiments performed with the parameters
$
K=128,\; L=2,$ and $h_t = 10^{-3}
$.
Figure~\ref{fig:init} shows the initial profiles of the three components
\(u_1,u_2,u_3\).  
Figure~\ref{fig:time} displays the corresponding profiles at several
increasing physical times \(t=10,50,100\), illustrating how the solution approaches a spatially
homogeneous equilibrium state.
\begin{figure}[ht]
	\includegraphics[scale=0.45]{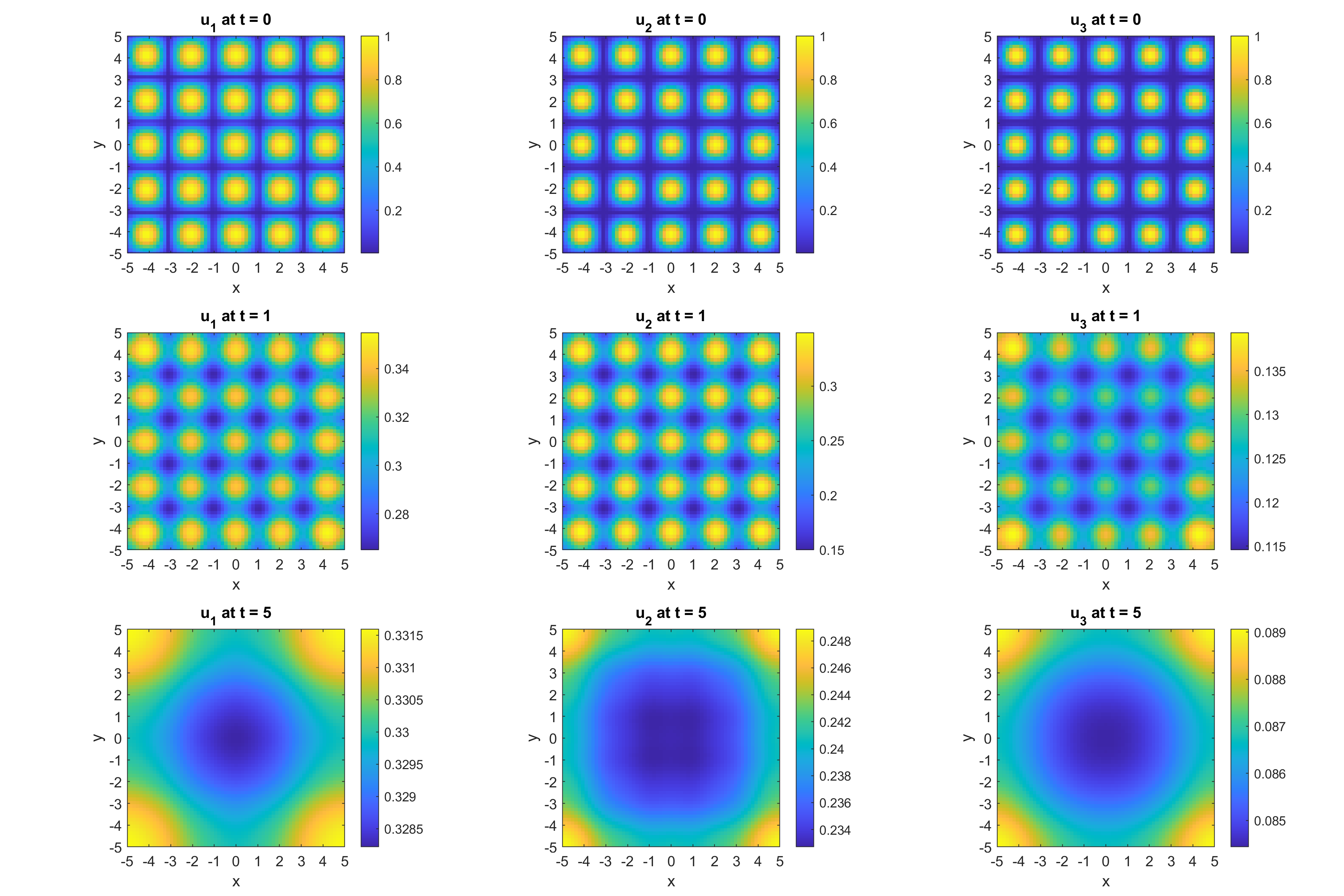}
	\centering
	\caption{
		Mid-plane slices (\(z=0\)) of the initial data and $(u_1,u_2,u_3)$ at times
		\(t = 1, 5\).
		%	Parameters: \(K=64\), \(L=2\), \(h_t=10^{-2}\).
	}
	\label{fig:init}
\end{figure}

\begin{figure}[ht]
	\includegraphics[scale=0.45]{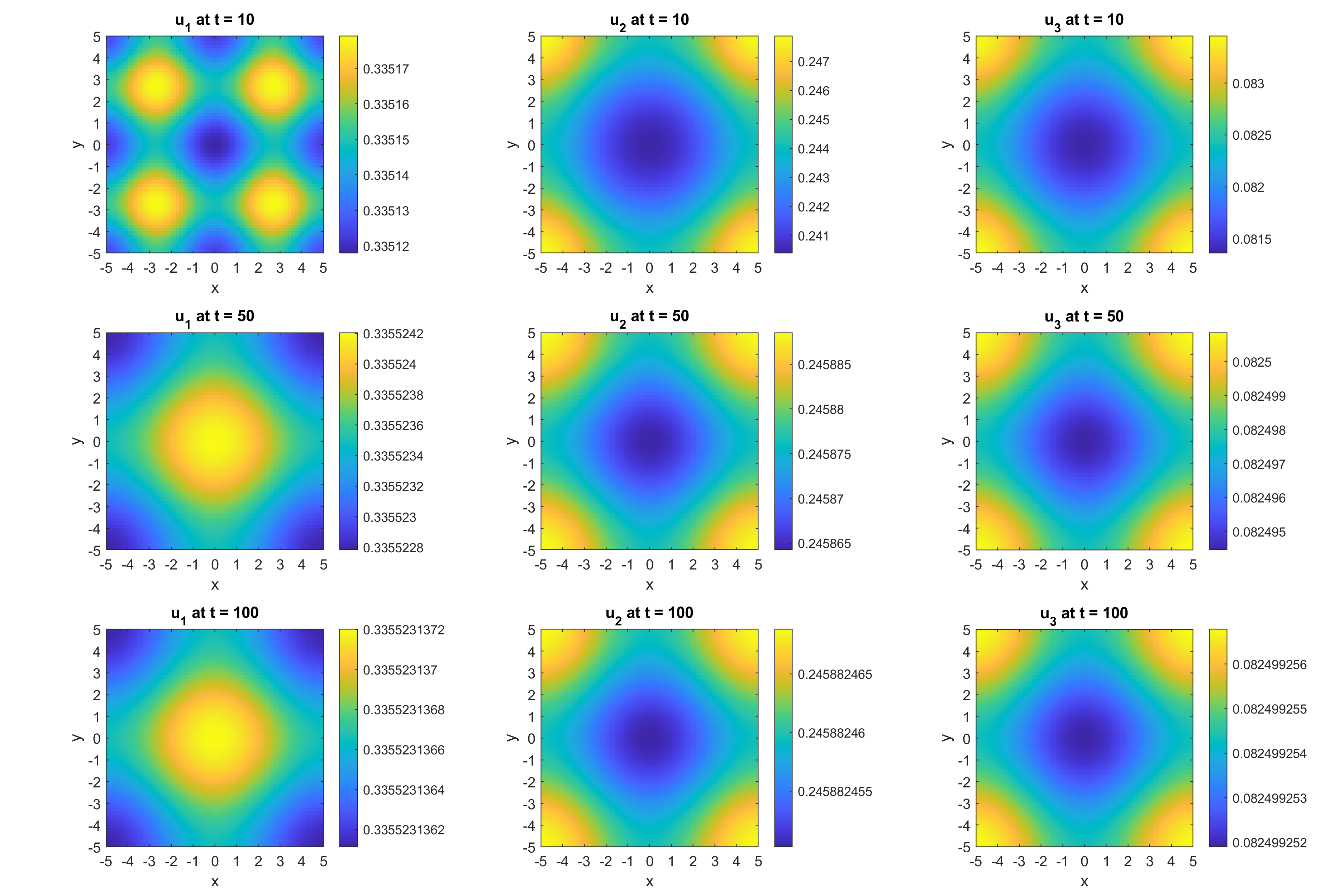}
		\centering
	\caption{
		Mid-plane slices (\(z=0\)) of $(u_1,u_2,u_3)$ at times
		\(t = 10, 50, 100\).
		%		Parameters: \(K=64\), \(L=2\), \(h_t=10^{-2}\).
	}
	\label{fig:time}
\end{figure}

\medskip
\noindent
Moreover, we compute the numerical solution up to a large final time and compare it to the equilibrium $(u_{1,\infty},u_{2,\infty},u_{3,\infty})$ {computed as in the classical case}.  
More precisely, we monitor the quantities
\[
E_i(t)
=\frac{\bigl\|u_i(t,\cdot)-u_{i,\infty}\bigr\|_{L^1(\Omega)}}{|u_{i,\infty}||\Omega|},
\qquad i=1,2,3.
\]
As illustrated in Figure~\ref{fig:L1cv}, the errors $E_i(t)$ decay rapidly to zero, showing that the three components converge numerically to the equilibrium values.

\begin{figure}[ht]
	\centering
	\includegraphics[scale=0.55]{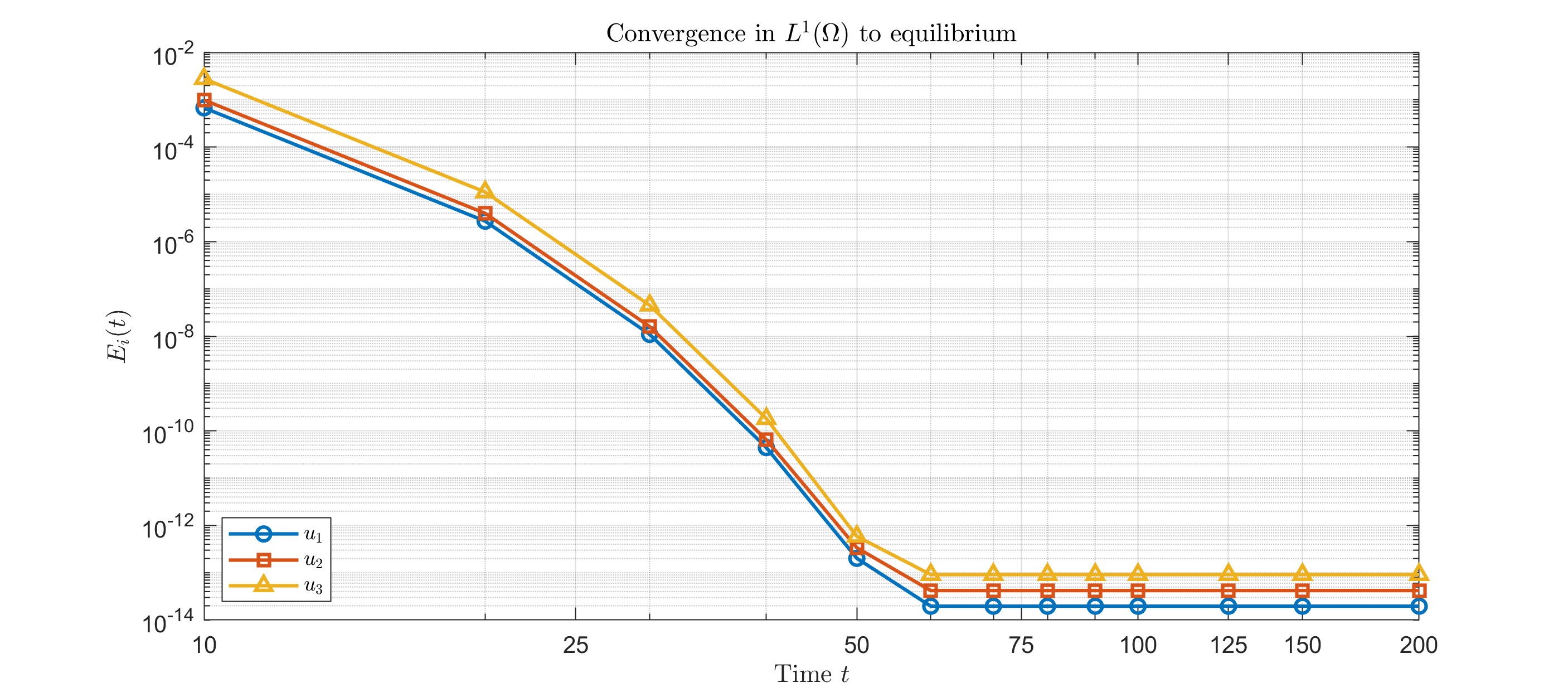}
	\caption{
		Relative $L^1(\Omega)$-errors $E_i(t)$ for $u_1$, $u_2$ and $u_3$ over time.
	}
	\label{fig:L1cv}
\end{figure}

Finally, we compute the \(L^2 \)-norm of \( u_1 \), \( u_2 \)  and \( u_3 \) at \( t_{\text{final}} = 5 \times 10^{5} \) for different values of \((K,L,h_t)\). 
The computed norms are presented in Table \ref{tab:l2norms}. 
As shown, the norms remain relatively constant across different values of  \((K,L,h_t)\), suggesting that the numerical solution is consistent.


\begin{table}
		\centering
	\begin{tabular}{ccc r r r}
		\hline\noalign{\smallskip}
		$K$ & $L$ & $h_t$ &
		$\|u_1\|_{L^2(\Omega)}$ & $\|u_2\|_{L^2(\Omega)}$ & $\|u_3\|_{L^2(\Omega)}$ \\
		\noalign{\smallskip}\hline\noalign{\smallskip}
		64  & 2 & $10^{-2}$ & 10.61017457 & 7.775486183 & 2.608856318 \\
		128 & 2 & $10^{-3}$ & 10.61017512 & 7.775486206 & 2.608856294 \\
		256 & 2 & $10^{-3}$ & 10.61017513 & 7.775486199 & 2.608856367 \\
		384 & 3 & $10^{-2}$ & 10.61017552 & 7.775486160 & 2.608856380 \\
		512 & 4 & $10^{-3}$ & 10.61017567 & 7.775486144 & 2.608856395 \\
		\noalign{\smallskip}\hline
	\end{tabular}
		\caption{$L^2$-norms of $(u_1,u_2,u_3)$ at $t_{\text{final}}=5\times 10^{5}$.}
	\label{tab:l2norms}
\end{table}

\bigbreak
In \textit{summary}, we have numerically investigated the global existence of
nonnegative solutions to System~\eqref{System-FR-Num}  in regimes not covered by current analytical tools. Our numerical simulations suggest that the system admits nonnegative
solutions that can be accurately computed over very large times. In addition, the numerical results confirm
the convergence of all components to the expected equilibrium state. However, the rigorous question of whether such solutions exist for all times is still open, and we intend to further explore this theoretical aspect in future work.

%\begin{acknowledgements}
%\end{acknowledgements}



\end{document}